\documentclass[times]{article}

\usepackage{hyperref}
\hypersetup{
	colorlinks=true,
	linkcolor=blue,
	citecolor=blue}
\usepackage{mathrsfs}
\usepackage{tikz}
\usepackage{float}
\usepackage{amsfonts,amssymb,amsbsy,amsmath,amsthm}
\usepackage[margin=1in]{geometry}
\usepackage{booktabs}
\usepackage{authblk}
\usepackage{titling}
\usepackage{enumitem}

\setlist[itemize]{topsep=7pt, itemsep=2pt, parsep=0pt, partopsep=0pt}

\def\dsqcup{\sqcup\mathchoice{\mkern-7mu}{\mkern-7mu}{\mkern-3.2mu}{\mkern-3.8mu}\sqcup}

\def\textfontii{\the\textfont\tw@}
\def\AmSTeX{{\textfontii A\kern-.1667em%
  \lower.5ex\hbox{M}\kern-.125emS}-\TeX\spacefactor1000 }
  
\long\def\@makecaption#1#2{%
  \vskip\abovecaptionskip
  \sbox\@tempboxa{\textbf{#1.}\hspace*{1em}#2}%
  \ifdim \wd\@tempboxa >\hsize
    \textbf{#1.}\hspace*{1em}#2\par
  \else
    \global \@minipagefalse
    \hb@xt@\hsize{\hfil\box\@tempboxa\hfil}%
  \fi
  \vskip\belowcaptionskip}

\newtheoremstyle{nonitalic}
  {\topsep}       
  {\topsep}       
  {\upshape}      
  {}              
  {\bfseries}     
  {.}             
  {.5em}          
  {}              

\theoremstyle{nonitalic}

\newtheorem{theorem}{Theorem}[section]

\newtheorem{corollary}[theorem]{Corollary}
\newtheorem{proposition}[theorem]{Proposition}
\newtheorem{definition}[theorem]{Definition}
\newtheorem{example}[theorem]{Example}

\newcommand\acks{\section*{Acknowledgements}}

\newtheorem{remark}[theorem]{Remark}

\newcommand{\sxl}[1]{\textcolor{black}{#1}}
\newcommand{\fss}[1]{\textcolor{black}{#1}}
\newcommand{\svwt}[1]{\textcolor{black}{#1}}
\newcommand{\fa}[1]{\textcolor{black}{#1}}
\newcommand{\fathree}[1]{\textcolor{black}{#1}}
\newcommand{\faa}[1]{\textcolor{black}{#1}}

\newcommand{\GX}{\mathscr{X}}
\newcommand{\GY}{\mathscr{Y}}

\newcommand{\comp}{\mathrm{\text{-}comp}}
\newcommand{\pa}{\mathrm{\text{-}par}}

\newcommand{\NCQSym}{\mathrm{NCQSym}}
\newcommand{\NCSym}{\mathrm{NCSym}}
\newcommand{\QSym}{\mathrm{QSym}}
\newcommand{\Sym}{\mathrm{Sym}}
\newcommand{\bM}{\mathbf{M}}
\newcommand{\bF}{\mathbf{F}}
\newcommand{\bx}{{\bf x}}
\newcommand{\by}{{\bf y}}
\newcommand{\bG}{{\mathbf{G}}}
\newcommand{\SSym}{{\mathfrak{S}{\rm Sym}}}
\newcommand{\arx}{{\begin{tikzpicture}
\draw[thick](0,0) circle (0.23);
\draw[thick,dashed,->](-0.23,0)--(0.23,0);
\end{tikzpicture}}
}

\newcommand{\ary}{{\begin{tikzpicture}
\draw[thick](0,0) circle (0.23);
\draw[thick,->](-0.23,0)--(0.23,0);
\end{tikzpicture}}
}
\newcommand{\arz}{{\begin{tikzpicture}
\draw[thick](0,0) circle (0.23);
\draw[thick,double,->](-0.23,0)--(0.23,0);
\end{tikzpicture}}}

\renewcommand{\maketitle}{
  \vspace*{4cm} 
  \begin{center}
    \Large\bfseries\thetitle
  \end{center}
  \vspace{1em}
  \begin{center}
    \large\theauthor
  \end{center}
}

\begin{document}

\title{\bf  Generalized Chromatic Functions}
\author{
Farid Aliniaeifard 
Shu Xiao Li, and
Stephanie van Willigenburg
}
\date{}

\maketitle

\begin{abstract}
We define vertex-colourings for edge-partitioned digraphs, which unify the theory 
of $P$-partitions and proper vertex-colourings of graphs. We use our vertex-colourings 
to define generalized chromatic functions, which merge the chromatic symmetric and
 quasisymmetric functions of graphs and generating functions of $P$-partitions. Moreover, 
 numerous classical bases of symmetric and quasisymmetric functions, both in commuting
  and noncommuting variables, can be realized as special cases of our generalized chromatic
   functions. We also establish product and coproduct formulas for our functions. Additionally,
    we construct the new Hopf algebra of $r$-quasisymmetric functions in noncommuting
     variables, and apply our functions to confirm its Hopf structure, and establish
      natural bases for it.
\end{abstract}


\section{Introduction} 
\fa{The story of the theory of $P$-partitions started with MacMahon's work \cite{Mac} on plane partitions
at the start of the 20th century, and 60 years later  Stanley, in his Ph.D. thesis \cite{Sta71}, extended the notion of plane partitions to $P$-partitions (for a more complete history, see \cite[pp. 169–188]{Ge15}).} The problem that  MacMahon considered in his work on plane partitions was the same as counting the number of fillings of Young diagrams with nonnegative integers with a given sum such that the entries are weakly decreasing in each row and column. 
\begin{center}
\begin{tikzpicture} 
\node at (0,0){$4$};
\node at (0.5,0){$3$};
\node at (1,0){$3$};
\node at (1.5,0){$1$};

\node at (0,-0.5){$4$};
\node at (0.5,-0.5){$3$};
\node at (1,-0.5){$2$};

\node at (0,-1){$2$};
\node at (0.5,-1){$1$};
\end{tikzpicture} 
\end{center}
Stanley, in his work on the theory of $P$-partitions \cite{Sta71}, generalized MacMahon's idea and replaced Young diagrams with posets and the weakly decreasing relation with weakly and strictly decreasing relations. In this paper, we generalize Stanley's $P$-partitions to certain vertex-colourings of digraphs whose edges are coloured with three colours. Roughly speaking, we replace posets in the theory of $P$-partitions with digraphs, and in addition to the weakly and strictly decreasing relations, we have another relation related to the proper colourings of digraphs. More precisely,  we define
a \emph{proper} vertex-colouring of a digraph whose edges are coloured by three colours, identified by $\dasharrow,\rightarrow,$ and  $\Rightarrow$, to be a function $\kappa:V(G)\rightarrow \mathbb{P}$, where $\mathbb{P}$ is the set of positive integers such that:
\begin{itemize} 
\item[(i)] If $a\dasharrow b$ in $G$, then $\kappa(a)\neq \kappa(b)$.
\item[(ii)] If $a\rightarrow b$ in $G$, then $\kappa(a)<\kappa(b)$.
\item[(iii)] If $a\Rightarrow b$ in $G$, then $\kappa(a)\leq \kappa(b)$.
\end{itemize} 
Let $\mathcal{C}(G)$ denote the set of all proper vertex-colourings of the edge-coloured digraph $G$. 
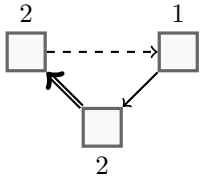
\begin{figure}[H]
\centering
 \begin{tikzpicture}[
roundnode/.style={circle, draw=black!60, fill=gray!10, very thick, minimum size=7mm},
squarednode/.style={rectangle, draw=black!60, fill=gray!5, very thick, minimum size=5mm},]
 \node[squarednode, label=above:$2$](1) at (1,1){};
  \node[squarednode,  label=below:$2$](3) at (2,0){};
    \node[squarednode,  label=above:$1$](2) at (3,1){};
    \draw[thick, dashed,->](1)--(2);
       \draw[ thick,->](2)--(3);
        \draw[ thick,double,->](3)--(1);
 \end{tikzpicture}
\caption{A proper vertex-colouring of an edge-coloured digraph}
\label{ex1}
\end{figure}
In the above definition, by replacing the infinite paintbox of colours  $\mathbb{P}$ with a finite paintbox of colours $[p]=\{1,2,\dots,p\}$, we generalize the classic proper vertex-colouring and weak/strong proper colouring in \cite{STW}.  The \emph{generalized chromatic number} of an edge-coloured digraph $G$ is the smallest number of colours we need to make a proper vertex-colouring of $G$.  Let $P(G,p)$ denote the number of ways that one can properly colour the vertices of  $G$ with $p$ colours. We show that $P(G,p)$ is a polynomial in $p$ (see Theorem \ref{thm:gcp}) and call it the \emph{generalized chromatic polynomial} of $G$.

 Stanley in \cite{Sta} defined the chromatic symmetric function of a finite simple graph. This symmetric function gives Birkhoff's chromatic polynomial by setting the first \faa{$p$} variables to $1$ and all others to $0$. There are two main conjectures regarding chromatic symmetric functions that have been open for more than 25 years---the {Tree Conjecture} \cite{Sta} and the $(3+1)$-free Conjecture \cite{SS}. Later on, refining the conjectures and seeking classical properties, other chromatic symmetric and quasisymmetric functions emerged, such as the chromatic quasisymmetric functions of graphs and digraphs \cite{E17, SW}, the extended chromatic symmetric function of a graph \cite{CS}, and the $k$-balanced chromatic quasisymmetric function of a graph  \cite{H}.

Let $\mathbb{Q}[[x_1,x_2,\dots]]$ be the set of all power series in commuting variables \fathree{$x=\{ x_1,x_2, \dots\}$.} We define the \emph{generalized chromatic function} of an edge-coloured digraph $G$ to be the bounded degree power series 
$$\GX_G(x)=\sum_{\kappa\in \mathcal{C}(G)} x_\kappa,$$where $x_\kappa=\prod_{a\in V(G)}x_{\kappa(a)}$. Note that the generalized chromatic polynomial $P(G,p)$ is equal to the generalized chromatic function of $G$ \fathree{after} setting the first $p$ variables to $1$ and all others to $0$. \sxl{In the theory of $P$-partitions, a labelled poset} $P$ corresponds to a quasisymmetric function $F_{P}$ \cite{Sta1}. Since each $P$-partition is a proper vertex-colouring of an edge-coloured digraph (see Section \ref{proper colourings}), we can \sxl{view $F_{P}$ as} a generalized chromatic function. Moreover, Stanley's chromatic symmetric functions \cite{Sta}, extended chromatic symmetric functions \cite{CS}, and chromatic quasisymmetric functions in \cite{E17,H, SW} are the generalized chromatic functions of certain edge-coloured digraphs (see Section \ref{others}). Therefore, we can merge all of these chromatic functions and the generating functions of $P$-partitions, $F_{P}$, into one object.

As we will see, generalized chromatic functions are quasisymmetric functions. There are many well-known bases for the Hopf algebras of symmetric and quasisymmetric functions. \fathree{Any of these bases can be expressed as a family of generalized chromatic functions. This offers a significant advantage because if one finds a product, coproduct, antipode, etc., formula for generalized chromatic functions, then this impacts knowledge about all these bases of the Hopf algebras of symmetric and quasisymmetric functions. For example, we give a generic coproduct formula for generalized chromatic functions.}

We also study a chain of Hopf algebras starting from symmetric functions and ending with quasisymmetric functions. We present different bases for the Hopf algebras in the chain using generalized chromatic functions. 

We then change gears to the noncommuting world and define \fa{$\GY_{\bG}(\bx)$}, the generalized chromatic function of a labelled edge-coloured digraph \fa{$\bG$} in noncommuting variables \fathree{$\bx=\{\bx_1,\bx_2,\dots\}$}. We extend most of the results in commuting variables to noncommuting variables.

\fathree{Moreover, we naturally expand the fundamental basis of the Malvenuto-Reutenauer Hopf algebra to the fundamental basis of the Hopf algebra of quasisymmetric functions in noncommuting variables. }

The Hopf algebra of $r$-quasisymmetric functions was defined in \cite{Hi} by Hivert. In \cite{GW}, Garsia and Wallach showed that the algebra of $r$-quasisymmetric functions is free over symmetric functions. In the last section, we introduce the Hopf algebra of $r$-quasisymmetric functions in noncommuting variables.

More precisely, our paper is structured as follows. In Section \ref{Hopf algebras}, we recall the background of symmetric, quasisymmetric, and $r$-quasisymmetric functions. In Section \ref{proper colourings}, we present some basic definitions in graph theory and then define some operators between edge-coloured digraphs. The vertex-colouring of an edge-coloured digraph is defined in Definition \ref{colouring}. We conclude Section  \ref{proper colourings} by showing that every $P$-partition is the proper vertex-colouring of an edge-coloured digraph in Proposition \ref{P-partition}.  In Section \ref{sec:GCF}, we introduce generalized chromatic functions in Definition \ref{GCF} and show that they are quasisymmetric functions, and then in Section \ref{others}, show that other chromatic symmetric and quasisymmetric functions are special cases of generalized chromatic functions. In Section \ref{prod-coprod}, the product and coproduct formulas for generalized chromatic functions are presented in Proposition \ref{prod} and Theorem \ref{coprod}. In Section \ref{bases}, we show that many well-known bases of symmetric, quasisymmetric, and $r$-quasisymmetric functions can be realized as special cases of generalized chromatic functions of edge-coloured digraphs. Moreover, in Theorem \ref{thm:gcp}, we show that the generalized chromatic polynomial of an edge-coloured digraph is indeed a polynomial. In Section  \ref{Hopf algebras noncomm}, we recall the background of symmetric and quasisymmetric functions in noncommuting variables and then define $r$-quasisymmetric functions in noncommuting variables. In Section \ref{labelled}, we introduce generalized chromatic symmetric functions in noncommuting variables in Definition \ref{GCFn}, and the product and coproduct formulas for them are presented in Propositions \ref{prodn} and \ref{coprodn}.  In Sections \ref{basesn} and \ref{MR}, we show that several bases for symmetric functions in noncommuting variables are the symmetrizations of certain generalized chromatic functions and give several bases for quasisymmetric functions in noncommuting variables, including its fundamental basis, which contains the fundamental basis of the Malvenuto-Reutenauer Hopf algebra. We conclude by showing that the set of $r$-quasisymmetric functions in noncommuting variables is a Hopf algebra in Theorem \ref{ncqsymr} and constructing the $r$-dominant monomial and upper-fundamental bases of the Hopf algebra of $r$-quasisymmetric functions in noncommuting variables in Proposition \ref{bncqsymr}.

\section{{Symmetric Functions and Generalizations}}\label{Hopf algebras}
This section introduces the Hopf algebras of symmetric, quasisymmetric, and $~r$-quasisymmetric functions. The bases of these Hopf algebras are indexed by partitions, compositions, and $r$-compositions, respectively.  We begin by recalling the definitions and notation related to these combinatorial objects.

\subsection{Partitions, compositions, and $r$-compositions}

A \emph{composition} $\alpha$ of $n$, denoted $\alpha\vDash n$, is a list of positive integers whose sum is $n$. 
 Given a composition $\alpha=(\alpha_1,\alpha_2,\ldots, \alpha_k),$ each $\alpha_i$ is called a \emph{part} of $\alpha$, the \emph{size} of $\alpha$ is $|\alpha|=\alpha_1+\alpha_2+\cdots+\alpha_k$, and the \emph{length} of $\alpha$ is $k$. For convenience we denote by $\emptyset$ the unique composition of size and length zero.
 If $\alpha=(\alpha_1,\alpha_2,\ldots,\alpha_k)\vDash n$, then we define 
 $$
 {\rm set}(\alpha)=\{\alpha_1,\alpha_1+\alpha_2,\dots,\alpha_1+\dots+\alpha_{k-1} \}\subseteq [n-1].
 $$
  For example, $(2,1,2)$ is a composition of $5$ with length $3$ and ${\rm set}(\alpha)=\{2,3\}$. 
For compositions $\alpha$ and $\beta$ of $n$, we write $\alpha\leq \beta$ and say $\alpha$ \emph{coarsens} $\beta$ (or $\beta$ \emph{refines} $\alpha$) if ${\rm set}(\alpha)\subseteq {\rm set}(\beta)$.

A \emph{partition} $\lambda=(\lambda_1,\lambda_2,\dots,\lambda_k)$ of $n$, denoted $\lambda \vdash n$, is a weakly decreasing composition. 
Let  $m_i$ be the number of parts of $\lambda$ that are equal to $i$. Let $\lambda^!=m_1!m_2!\cdots m_n!$, and let $\lambda!=\lambda_1!\lambda_2!\cdots \lambda_k!$.  We sometimes write $$\lambda=(n^{m_n}, (n-1)^{m_{n-1}}, \ldots,1^{m_1}).$$ 
\fa{The \emph{dominance order} of partitions is defined as follows.} For partitions $\lambda=(\lambda_1,\lambda_2,\ldots,\lambda_k)$ and $\mu=(\mu_1,\mu_2,\ldots,\mu_l)$ of $n$, we write $\mu\preceq \lambda$  if $k\leq l$ and for every $1\leq i \leq k$, $$ \mu_1+\mu_2+\cdots+\mu_i \leq  \lambda_1+\lambda_2+\cdots+\lambda_i .$$
For example, $(3,1,1,1) \preceq (3,2,1)$.

Let $r$ be a positive integer or infinity. A pair $(\beta,\mu)$ is called an \emph{$r$-composition} of $|\beta|+|\mu|$ when:
\begin{itemize} 
\item[(i)] $\beta$ is a composition whose parts are at least $r$.
\item[(ii)] $\mu$ is a partition whose parts are strictly smaller than $r$.
\end{itemize}  For example, $((7,4,5),(3,2))$ is a $4$-composition of $21$.

\subsection{Quasisymmetric functions}
The Hopf algebra of quasisymmetric functions was formally introduced by Gessel \cite{Ge} in 1984. From this concept, a whole research area emerged;
a history can be found in \cite[Introduction]{LMvW}. 

Recall that $\mathbb{Q}[[x_1,x_2,\ldots]]$ is the algebra of formal power series  in infinitely many commuting variables \fathree{$x=\{x_1,x_2,\ldots\}$} over $\mathbb{Q}$. Let $\mathfrak{S}_n$ be the group of all permutations of $[n]$. Let $\mathfrak{S}_\infty=\sqcup_{n\geq 0}\mathfrak{S}_n$. We identify a permutation $\sigma\in \mathfrak{S}_n\subseteq \mathfrak{S}_\infty$ with a bijection of the positive integers by defining $\sigma(i)=i$ if $i>n$. 

\begin{definition}
A \emph{quasisymmetric function} is a formal power series $f\in \mathbb{Q}[[x_1,x_2,\ldots]]$ such that:
\begin{itemize}
\item[(i)] The degrees of the monomials in $f$ are bounded.
\item[(ii)] For every composition $(\alpha_1,\alpha_2,\ldots,\alpha_k)$, all monomials $x_{i_1}^{\alpha_1} x_{i_2}^{\alpha_2}  \cdots x_{i_k}^{\alpha_k}$ in $f$ with indices $i_1<i_2<\cdots < i_k$ have the same coefficient.
\end{itemize}
The set of all quasisymmetric functions is denoted by $\QSym(x)$. 
\end{definition}

The vector space $\QSym(x)$ is a Hopf algebra, where its product is the same as the product of the formal power series and its coproduct  $\Delta$ is  defined as follows (for more details see \cite[p. 142]{GR}). Consider the linear order on two sets of commuting variables $(x,y)=(x_1<x_2<\cdots < y_1<y_2<\cdots)$, and inject $\QSym(x)\otimes \QSym(y)$ into $\mathbb{Q}[[x,y]]$ by identifying every $f\otimes g \in \QSym(x)\otimes \QSym(y)$ with $fg\in \mathbb{Q}[[x,y]]$. We then have that
$$\QSym(x,y)\subseteq \QSym(x)\otimes \QSym(y).$$ We can define $\Delta:\QSym(x) \rightarrow \QSym(x)\otimes \QSym(x)$ as the composite of the following maps. 
$$
\begin{array}{cccccc}
\QSym(x)  \cong & \QSym(x,y) \rightarrow &\QSym(x)\otimes \QSym(y) \cong& \QSym(x)\otimes \QSym(x)\\
f \mapsto & f(x_1,x_2,\ldots,y_1,y_2,\ldots)  & & 
\end{array}
$$

 For more details about the Hopf algebra of quasisymmetric functions and its well-known bases, see \cite{Sta1}.

\subsection{Symmetric functions}

 The Hopf algebra of symmetric functions is a Hopf subalgebra of $\QSym(x)$. 
\begin{definition}
A \emph{symmetric function}  is a formal power series $f\in \mathbb{Q}[[ x_1,x_2,\ldots ]]$ such that:
\begin{itemize}
\item[(i)] The degrees of the monomials in $f$ are bounded.
\item[(ii)] For any permutation $\sigma\in\mathfrak{S}_\infty$, $$\sigma \cdot f(x_1,x_2,\ldots)=f(x_{\sigma(1)},x_{\sigma(2)}, \ldots)=f(x_1,x_2,\ldots).$$
\end{itemize}
 The set of all symmetric functions is  denoted by $\Sym(x)$. 
 \end{definition}
 For more details about the Hopf algebra of symmetric functions and its well-known bases, see \cite{Sta1}.

\subsection{$r$-quasisymmetric functions}

Hivert introduced the Hopf algebra of $r$-quasisymmetric functions in  \cite{Hi}, which is a Hopf subalgebra of $\QSym(x)$.  

For each $r$-composition $(\beta,\mu)$ where $\beta=(\beta_1,\beta_2,\ldots,\beta_k)$ and $\mu=(\mu_1,\mu_2,\ldots,\mu_l)$, define the \emph{$r$-dominant monomial function} to be
$$M_{(\beta,\mu)}=\sum x_{i_1}^{\beta_1}x_{i_2}^{\beta_2}\cdots x_{i_k}^{\beta_k} x_{i_{k+1}}^{\mu_1}x_{i_{k+2}}^{\mu_2} \cdots x_{i_{k+l}}^{\mu_l},$$
where the sum is over all distinct positive integers $i_1, i_2, \ldots, i_{k+l}$ such that \fathree{only the first $k$ indices are required to be in strictly increasing order, that is} $i_1<i_2<\cdots < i_k$.
Define
$$
{\rm QSym}^r(x)=\bigoplus_{n\geq 0}{\rm QSym}_n^r(x),
$$ 
where 
$${\rm QSym}_n^r(x)=\mathbb{Q}\text{-span}\{ M_{(\beta,\mu)}: (\beta,\mu) \text{~ is an $r$-composition of $n$}\}.$$
We have that 
$${\rm QSym}(x)={\rm QSym}^1(x)  \supset {\rm QSym}^2(x)  \supset \cdots \supset {\rm QSym}^\infty(x)={\rm Sym}(x).$$
In Section \ref{bases},  we present different bases for the Hopf algebras in this chain.

\section{Proper Colourings}\label{proper colourings}

In graph theory, there are many families of graph colourings. These colourings are usually defined by setting specific constraints on the colours of vertices or edges of a graph. In this paper, we consider particular vertex-colourings of a digraph to unify several combinatorial constructions. In our colouring, the constraints on the colours of the vertices are subject to an edge-colouring of the digraph. We recall the definitions and notation in graph theory that we need throughout this paper.

A \emph{simple digraph} $G=(V(G),E(G))$ is a digraph with no loops, and for any distinct vertices $a$ and $b$, there can be at most one directed edge
from $a$ to $b$. \fathree{Note that for distinct vertices  $a$ and $b$ of a simple digraph, we can have two edges one from $a$ to $b$ and one from $b$ to $a$.}
Throughout this paper, all digraphs are finite and simple, and we usually use $G$ to denote a simple digraph.  \fa{The \emph{underlying graph} of $G$ is the undirected graph whose vertex set is the vertex set of $G$ and its edge set is $\{ab: \text{$(a,b)$ or $(b,a)$ is an edge of $G$}\}$}. A \emph{directed cycle} is a digraph $G$ with the vertex set $V(G)=\{a_1,a_2,\ldots, a_n\}$ and the edge set $E(G)=\{(a_i,a_{i+1}): i=1, 2, \ldots, n-1\}\fa{\uplus} \{(a_{n},a_1)\}$. A \emph{directed path} is a digraph $G$ with the vertex set $V(G)=\{a_1,a_2,\ldots, a_n\}$ and the edge set $E(G)=\{(a_i,a_{i+1}): i=1,2,\ldots, n-1\}$.  A \emph{complete digraph} is a digraph $G$ with the vertex set $V(G)=\{a_1,a_2,\ldots, a_n\}$ and the edge set $E(G)=\{(a_i,a_{j}): 1\leq i<j\leq n\}$. 

Let $S$ and $S'$ be subsets of $\mathbb{P}$, the set of positive integers. We regard the elements in $S$ and $S'$ as colours.  Let $G$ be a digraph. \fa{An} \emph{$S$-vertex-colouring} of $G$ is a function $\kappa$ that assigns a colour in $S$ to each vertex of the digraph $G$. By a vertex-colouring, without mentioning the set $S$, we mean a $\mathbb{P}$-vertex-colouring. 
An \emph{$S'$-edge-colouring} of $G$ is a function $\kappa'$ that assigns a colour in $S'$ to each edge of the digraph $G$. 
An \emph{$S'$-edge-coloured} digraph is a digraph $G$ together with an $S'$-edge-colouring $\kappa'$ of $G$.  {\emph{ Throughout this paper, we only consider $S'$-edge-colourings of digraphs where $|S'|=3$, so we only have three types of coloured edges in a digraph, which are denoted by $\dashrightarrow ,\rightarrow$, and $\Rightarrow$.
}}
\begin{figure}[H]
\begin{center}
  \begin{tikzpicture}[
roundnode/.style={circle, draw=black!60, fill=gray!10, very thick, minimum size=7mm},
squarednode/.style={rectangle, draw=black!60, fill=gray!5, very thick, minimum size=5mm},]
 \node[squarednode](1) at (1,1){};
  \node[squarednode](3) at (2,0){};
    \node[squarednode](2) at (3,1){};
    \draw[thick, dashed,->](1)--(2);
       \draw[ thick,->](2)--(3);
        \draw[ thick,double,->](3)--(1);
 \end{tikzpicture}
   \caption{An edge-coloured digraph}
   \label{ex1edge}
\end{center}
\end{figure}
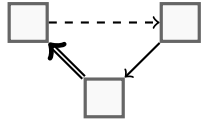

Consider edge-coloured digraphs $G_1$ with the edge-colouring $\kappa_1'$ and $G_2$ with the edge-colouring $\kappa_2'$. 
  The \emph{disjoint union} of the edge-coloured digraphs $G_1$ and $G_2$, denoted $G_1 \uplus G_2$,
 is an edge-coloured digraph $G$ together with an edge-colouring $\kappa'$ such that:
 \begin{itemize} 
\item[(i)] The vertex set of $G$ is the disjoint union of the vertex sets of $G_1$ and $G_2$.
\item[(ii)] The  edge set of $G$ is the disjoint union of the edge sets of $G_1$ and $G_2$. 
\item[(iii)]  $\kappa'(e)=\kappa'_j(e)$ if $e\in E(G_j)$ for $j=1,2$. 
 \end{itemize}
 
 \begin{figure}[H]
 \begin{center}
  \begin{tikzpicture}[
roundnode/.style={circle, draw=black!60, fill=gray!10, very thick, minimum size=7mm},
squarednode/.style={rectangle, draw=black!60, fill=gray!5, very thick, minimum size=5mm},]
 \node[squarednode](1) at (1,1){};
  \node[squarednode](3) at (2,0){};
    \node[squarednode](2) at (3,1){};
    \draw[thick, dashed,->](1)--(2);
       \draw[ thick,->](2)--(3);
        \draw[ thick,double,->](3)--(1);
        \node at (4,0.5){$\biguplus$};
 \end{tikzpicture}~~  \begin{tikzpicture}[
roundnode/.style={circle, draw=black!60, fill=gray!10, very thick, minimum size=7mm},
squarednode/.style={rectangle, draw=black!60, fill=gray!5, very thick, minimum size=5mm},]
 \node[squarednode](1) at (0,1.5){};
  \node[squarednode](2) at (1.5,1.5){};
    \node[squarednode](3) at (1.5,0){};
      \node[squarednode](4) at (0,0){};
    \draw[ thick,double ,->](1)--(2);
     \draw[ thick,dashed ,->](2)--(3);
      \draw[ thick,dashed ,->](3)--(4);
       \draw[thick,solid ,->](4)--(1);
      \node at (2.5,0.5){$=$};
 \end{tikzpicture} ~~
 \begin{tikzpicture}[
roundnode/.style={circle, draw=black!60, fill=gray!10, very thick, minimum size=7mm},
squarednode/.style={rectangle, draw=black!60, fill=gray!5, very thick, minimum size=5mm},]
 \node[squarednode](1) at (1,1){};
  \node[squarednode](3) at (2,0){};
    \node[squarednode](2) at (3,1){};
    \draw[thick, dashed,->](1)--(2);
       \draw[ thick,->](2)--(3);
        \draw[ thick,double,->](3)--(1);
 \end{tikzpicture}  \begin{tikzpicture}[
roundnode/.style={circle, draw=black!60, fill=gray!10, very thick, minimum size=7mm},
squarednode/.style={rectangle, draw=black!60, fill=gray!5, very thick, minimum size=5mm},]
 \node[squarednode](1) at (0,1.5){};
  \node[squarednode](2) at (1.5,1.5){};
    \node[squarednode](3) at (1.5,0){};
      \node[squarednode](4) at (0,0){};
    \draw[thick,double ,->](1)--(2);
     \draw[thick,dashed ,->](2)--(3);
      \draw[thick,dashed ,->](3)--(4);
       \draw[thick,,solid ,->](4)--(1);
 \end{tikzpicture}
 \caption{The disjoint union of two edge-coloured digraphs}
 \label{dj}
\end{center}
\end{figure}
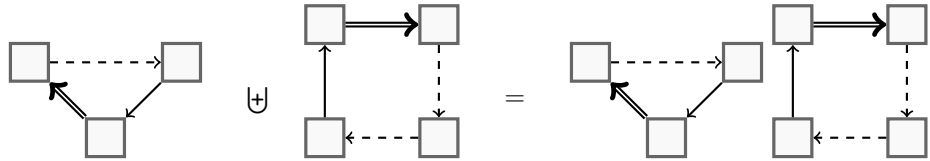
Also, the \emph{dashed} (\emph{solid} and \emph{double}, respectively) \emph{sum} of the edge-coloured digraphs $G_1$ and $G_2$ is an edge-coloured digraph $G$ together with an edge-colouring $\kappa'$ such that:
 \begin{itemize}
 \item[(i)] The  vertex set of $G$ is the disjoint union of the vertex sets of $G_1$ and $G_2$.
 \item[(ii)] The edge set of $G$ is the disjoint union of the edge set of $G_1$, edge set of $G_2$, and $\{(a,b):a \text{~is a vertex of $G_1$ and $b$ is a vertex of}$ $G_2\}$.  
 \item[(iii)] For every edge $e=(a,b)$ in $E(G)$, $\kappa'(e)=\kappa'_j(e)$ if $e\in E(G_j)$ for $j=1,2,$ and if $a$ is a vertex of $G_1$ and $b$ is a vertex of $G_2$, then $e$ is a dashed (solid and double, respectively)  edge.
 \end{itemize} 
 The dashed, solid, and double sums of the edge-coloured digraphs $G_1$  and $G_2$ are denoted  by
  $$
  \begin{tikzpicture}[
roundnode/.style={circle, draw=black!60, fill=gray!10, very thick, minimum size=5mm},
squarednode/.style={rectangle, draw=black!60, fill=gray!5, very thick, minimum size=5mm},]
    \node(1) at (0.2,0){$G_1$};
      \node(2) at (1.3,0){$G_2,$};
        \node(a) at (0.7,0){$\arx$}; 
   \end{tikzpicture} 
   \begin{tikzpicture}[
roundnode/.style={circle, draw=black!60, fill=gray!10, very thick, minimum size=5mm},
squarednode/.style={rectangle, draw=black!60, fill=gray!5, very thick, minimum size=5mm},]
    \node(1) at (0.2,0){$G_1$};
      \node(2) at (1.3,0){$G_2,$};
         \node(a) at (0.7,0){$\ary$}; 
           \node (b) at (2.3,0){and}; 
   \end{tikzpicture} 
   \begin{tikzpicture}[
roundnode/.style={circle, draw=black!60, fill=gray!10, very thick, minimum size=5mm},
squarednode/.style={rectangle, draw=black!60, fill=gray!5, very thick, minimum size=5mm},]
    \node(1) at (0.2,0){$G_1$};
      \node(2) at (1.3,0){$G_2,$};
       \node(a) at (0.7,0){$\arz$}; 
   \end{tikzpicture}
$$ respectively.   

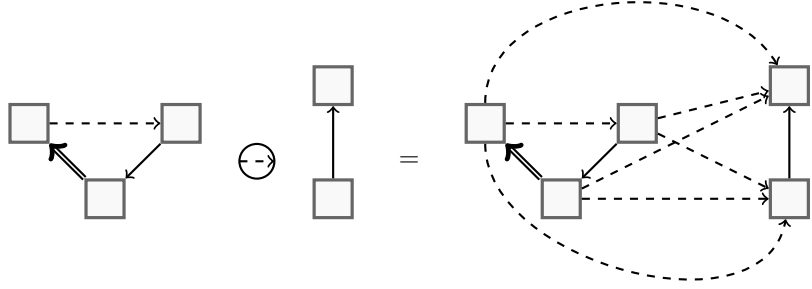
\begin{figure}[H]
\begin{center}
  \begin{tikzpicture}[
roundnode/.style={circle, draw=black!60, fill=gray!10, very thick, minimum size=7mm},
squarednode/.style={rectangle, draw=black!60, fill=gray!5, very thick, minimum size=5mm},]
 \node[squarednode](a1) at (1,1){};
  \node[squarednode](a3) at (2,0){};
    \node[squarednode](a2) at (3,1){};
    \draw[thick, dashed,->](a1)--(a2);
       \draw[ thick,->](a2)--(a3);
        \draw[ thick,double,->](a3)--(a1);
        \node at (4,0.5){$\arx$};
  \node[squarednode](b1) at (5,1.5){};
      \node[squarednode](b4) at (5,0){};
       \draw[thick,solid ,->](b4)--(b1);
      \node at (6,0.5){$=$};
 
 \node[squarednode](1) at (7,1){};
  \node[squarednode](3) at (8,0){};
    \node[squarednode](2) at (9,1){};
    \draw[thick, dashed,->](1)--(2);
       \draw[ thick,->](2)--(3);
        \draw[ thick,double,->](3)--(1);
        
        \node[squarednode](11) at (11,1.5){};
         \node[squarednode](14) at (11,0){};
         
          \draw[thick,solid ,->](14)--(11);
          
          \draw[thick,dashed ,->](1) to [out=90,in=120] (11);
          \draw[thick,dashed ,->](2)--(11);
          \draw[thick,dashed ,->](3)--(11);
           \draw[thick,dashed ,->](1) to [out=-90,in=-100]  (14);
          \draw[thick,dashed ,->](2)--(14);
          \draw[thick,dashed ,->](3)--(14);
 \end{tikzpicture}  
  \caption{The dashed sum of two edge-coloured digraphs}
  \label{dash}
\end{center}
\end{figure}

We frequently use the edge-coloured digraphs in the following table.

\begin{table}[H]
\caption{Useful edge-coloured digraphs}
\label{ecd}
\centering
\begin{tabsize}
\begin{tabular}{|c|c|}
 \toprule
 \text{Notation} & \text{Expression}  \\
\midrule
 $C_n$ &  \text{The directed cycle with $n$ vertices and double edges}\\
 \midrule
 $P_n$ & \text{The directed path with $n$ vertices and solid edges}\\
\midrule
 $Q_n$ &  \text{The directed path with $n$ vertices and double edges}\\
\midrule
 $K_n$ &  \text{The complete digraph with $n$ vertices and dashed edges}\\
 \bottomrule
 \end{tabular}
 \end{tabsize}
\end{table}

\begin{figure}[H]
 \begin{center}
   \begin{tikzpicture}[
roundnode/.style={circle, draw=black!60, fill=gray!10, very thick, minimum size=7mm},
squarednode/.style={rectangle, draw=black!60, fill=gray!5, very thick, minimum size=5mm},]
 \node[squarednode](1) at (1.3,0){};
  \node[squarednode](2) at (2.6,0){};
    \node[squarednode](3) at (3.9,0){};
     \node at (2.6,-0.7){$C_3$};
    \draw[thick,double,->](1)--(2);
       \draw[ thick,double,->](2)--(3);
        \draw[ thick,double,->](3)--(3.9,0.7)--(1.3,0.7)--(1);
 \end{tikzpicture} \quad 
   \begin{tikzpicture}[
roundnode/.style={circle, draw=black!60, fill=gray!10, very thick, minimum size=7mm},
squarednode/.style={rectangle, draw=black!60, fill=gray!5, very thick, minimum size=5mm},]
 \node[squarednode](1) at (1.3,0){};
  \node[squarednode](2) at (2.6,0){};
    \node[squarednode](3) at (3.9,0){};
     \node at (2.6,-0.7){$P_3$};
    \draw[thick,->](1)--(2);
       \draw[ thick,->](2)--(3);
 \end{tikzpicture} \quad
  \begin{tikzpicture}[
roundnode/.style={circle, draw=black!60, fill=gray!10, very thick, minimum size=7mm},
squarednode/.style={rectangle, draw=black!60, fill=gray!5, very thick, minimum size=5mm},]
 \node[squarednode](1) at (1.3,0){};
  \node[squarednode](2) at (2.6,0){};
    \node[squarednode](3) at (3.9,0){};
    \node at (2.6,-0.7){$Q_3$};
    \draw[thick,double,->](1)--(2);
       \draw[ thick,double,->](2)--(3);
 \end{tikzpicture} \quad
  \begin{tikzpicture}[
roundnode/.style={circle, draw=black!60, fill=gray!10, very thick, minimum size=7mm},
squarednode/.style={rectangle, draw=black!60, fill=gray!5, very thick, minimum size=5mm},]
 \node[squarednode](1) at (1.3,0){};
  \node[squarednode](2) at (2.6,0){};
    \node[squarednode](3) at (3.9,0){};
     \node at (2.6,-0.7){$K_3$};
    \draw[thick,dashed,->](1)--(2);
       \draw[ thick,dashed,->](2)--(3);
        \draw[ thick,dashed,->](1)--(1.3,0.7)--(3.9,0.7)--(3);
 \end{tikzpicture}
  \caption{The edge-coloured digraphs $C_3,P_3,Q_3, K_3$}
  \label{exgraphs}
  \end{center}
  \end{figure}
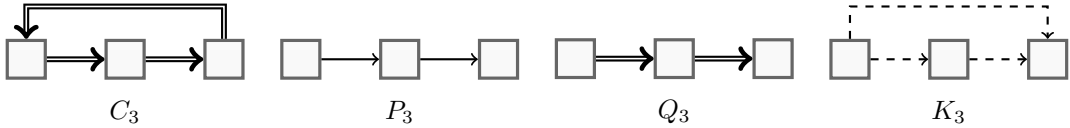

Let $(P,<_P)$ be a poset. The \emph{digraph associated  to $P$}, denoted $G_P$, is a digraph whose vertex set is the elements of $P$ and $(a,b)$ is an edge of $G_P$ if and only if $b$ \emph{covers} $a$ in $P$; that is, $a<_Pb$ and there is no $c\in P$ such that $a<_Pc<_Pb.$ \fathree{Note that the underlying graph of $G_P$ gives the Hasse diagram of $P$.} 

We now define our vertex-colouring of an edge-coloured digraph $G$, which plays an essential role in this paper. 
\begin{definition} \label{colouring}
A proper vertex-colouring of an edge-coloured digraph $G$ is a function $\kappa$ from $V(G)$ to $\mathbb{P}$ such that for any edge $(a,b)$ in $E(G)$:
\begin{itemize}
\item[(i)] If $a\dashrightarrow b$, then $\kappa(a)\neq \kappa(b)$.
\item[(ii)] If  $a \rightarrow b$, then $\kappa(a)< \kappa(b)$.
\item[(iii)] If  $a \Rightarrow b$, then $\kappa(a) \leq  \kappa(b)$.
\end{itemize} 
Let $\mathcal{C}(G)$ denote the set of all proper vertex-colourings of the edge-coloured digraph $G$. 
\end{definition}

\begin{remark}
 Note that for some edge-coloured digraphs $G$, $\mathcal{C}(G)$ is empty.
\end{remark}

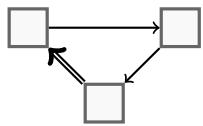
\begin{figure}[H]
\begin{center}
  \begin{tikzpicture}[
roundnode/.style={circle, draw=black!60, fill=gray!10, very thick, minimum size=7mm},
squarednode/.style={rectangle, draw=black!60, fill=gray!5, very thick, minimum size=5mm},]
 \node[squarednode](1) at (1,1){};
  \node[squarednode](3) at (2,0){};
    \node[squarednode](2) at (3,1){};
    \draw[thick,->](1)--(2);
       \draw[ thick,->](2)--(3);
        \draw[ thick,double,->](3)--(1);
 \end{tikzpicture}
   \caption{An edge-coloured digraph $G$ with $\mathcal{C}(G)=\emptyset$}
   \label{rem:emptycolouring}
\end{center}
\end{figure}

We next see that each $P$-partition corresponds to a proper vertex-colouring of $G_P$ where its edges are solid or double.

\paragraph{Proper colourings and $P$-partitions.}
 A \emph{labelled poset} is a partially ordered set $P$ whose underlying set is some finite subset of positive integers.   A \emph{$P$-partition} is a function $f: P\rightarrow \mathbb{P}$ such that for \fa{$a,b\in P$ with $a<_P b$}:
 \begin{itemize}
\item[(i)]  If   $a<_{\mathbb{Z}} b$, then 
 $f(a)\leq f(b)$.
 \item[(ii)] If  $a>_{\mathbb{Z}} b$, then 
 $f(a)< f(b)$. 
 \end{itemize}

The definitions of $P$-partitions and proper vertex-colourings of edge-coloured digraphs yield the following proposition. 
 
\begin{proposition}\label{P-partition}
Let $P$ be a poset, and let $G_P$ be the digraph associated to $P$. 
Then there is a bijection between the set of all $P$-partitions and the set of all proper colourings of the edge-coloured digraph $G$ where the digraph $G$ is isomorphic to $G_P$ and moreover: 
\begin{itemize}
\item[(i)] $a\Rightarrow b$ in $G$ if $a<_P b$ and $a<_{\mathbb{Z}} b$.
\item[(ii)]  $a\rightarrow b$ in $G$ if $a<_P b$ and $a>_{\mathbb{Z}} b$.
\end{itemize}
\end{proposition}

\section{Generalized Chromatic Functions}\label{sec:GCF}

\fa{
Throughout this paper, if $\kappa$ is any type of colouring of a graph, digraph, or edge-coloured digraph with vertex set $V$, we define 
\begin{equation}\label{eq:xkappa}
x_\kappa=\prod_{a\in V} x_{\kappa(a)}.
\end{equation} 
}
 We now introduce our main object of study. 

\begin{definition}\label{GCF} {(The generalized chromatic function of an edge-coloured digraph)}
 Let $G$ be an edge-coloured digraph. Define
$$
\GX_{G}(x,t)=\sum_{\kappa\in \mathcal{C}(G)}t^{{\rm asc}(\kappa)} x_\kappa,
$$
\fa{where $x_\kappa$ is defined by Equation \ref{eq:xkappa}} and 
$$ {\rm asc}(\kappa)=|\{(a,b)\in E(G):  \kappa(a)<\kappa(b)\}|.$$
The power series 
$$\GX_{G}(x)=\GX_{G}(x,1)$$ is called the  \emph{generalized chromatic function} of $G$.
\end{definition} 

\begin{remark}
Note that $\GX_{G}(x,t)=0$ when $ \mathcal{C}(G)=\emptyset$.
\end{remark}

If $G$ is an edge-coloured digraph and $a\dashrightarrow b$ is an edge of $G$, then the generalized chromatic function of $G$ is the sum of the generalized chromatic functions of two edge-coloured digraphs, one is obtained by deleting $a\dashrightarrow b$ in $G$ and replacing it with $a\rightarrow b$ and the other  is  obtained by deleting $a\dashrightarrow b$ in $G$ and replacing it with $b\rightarrow a$; that is
$$\GX_G(x)=\GX_{G-(a\dashrightarrow b)+(a\rightarrow b)}(x)+\GX_{G-(a\dashrightarrow b)+(b\rightarrow a)}(x).$$

 For vertices $a$ and $b$ of $G$, write $a\sim b$ if $a=b$ or there is a directed cycle $C$ with double edges in $G$ and $a,b\in V(C)$.  
 \fa{The} transitive closure of this relation gives an equivalence relation on the vertices of $G$. 
Let $[a_1], [a_2],\ldots,[a_s]$ be all equivalence classes. We have
$$
\GX_{G}(x)=\sum_{\kappa\in \mathcal{C}(G)} x_{\kappa(a_1)}^{|[a_1]|} x_{\kappa(a_2)}^{|[a_2]|} \cdots x_{\kappa(a_s)}^{|[a_s]|}.
$$
\fa{By the above equation} we can realize each generalized chromatic function as a sum of generating functions of weighted $P$-partitions \cite[Section 3]{AWvW}, and so the generalized chromatic function of an edge-coloured digraph is a quasisymmetric function, which we now state as a proposition.

\begin{proposition}\label{prop:gcfisqsym}
Let $G$ be an edge-coloured digraph. Then $\GX_G(x) \in \faa{\QSym(x)}$.
\end{proposition}

 For example, if $G$ is the following edge-coloured digraph
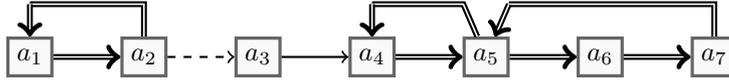
\begin{center} 
  \begin{tikzpicture}[
roundnode/.style={circle, draw=black!60, fill=gray!10, very thick, minimum size=7mm},
squarednode/.style={rectangle, draw=black!60, fill=gray!5, very thick, minimum size=5mm},]
 \node[squarednode](1) at (0,0){$a_1$};
  \node[squarednode](2) at (1.5,0){$a_2$};
    \node[squarednode](3) at (3,0){$a_3$};
      \node[squarednode](4) at (4.5,0){$a_4$};
        \node[squarednode](5) at (6,0){$a_5$};
          \node[squarednode](6) at (7.5,0){$a_6$};
            \node[squarednode](7) at (9,0){$a_7$};
     
\draw[thick,double,->](1)--(2);
\draw[thick,double](2)--(1.5,0.7)--(0,0.7);
\draw[thick,double,->](0,0.7)--(1);

\draw[thick,dashed,->](2)--(3);

\draw[thick,->](3)--(4);
\draw[thick,double,->](4)--(5);
\draw[thick,double](5)--(5.7,0.7)--(4.5,0.7);
\draw[thick,double,->](4.5,0.7)--(4);

\draw[thick,double,->](5)--(6);
\draw[thick,double,->](6)--(7);
\draw[thick,double](7)--(9,0.7)--(6.3,0.7);
\draw[thick,double,->](6.3,0.7)--(5);

  \end{tikzpicture}
    \end{center} 
then $[a_1]=\{a_1,a_2\}$, $[a_3]=\{a_3\}$, and $[a_4]=\{ a_4, a_5,a_6,a_7\}$ are  the equivalence classes, and
$$\GX_{G}(x)=\GX_{G-(a_2\dashrightarrow a_3)+(a_2 \rightarrow a_3)}(x)+\GX_{G-(a_2\dashrightarrow a_3)+(a_3 \rightarrow a_2)}(x).$$

\section{Other Chromatic Functions}\label{others}
Let $H=(V(H),E(H))$ be a finite simple graph. Throughout this paper, all graphs are finite and simple, and we usually use $H$ to denote a simple graph.  A \emph{proper} vertex-colouring of  $H$ is a function $\kappa$ from $V(H)$ to $\mathbb{P}$ such that if the vertices $a$ and $b$ are adjacent, then $\kappa(a)\neq \kappa(b)$. The set of all proper vertex-colourings of $H$ is denoted by $\mathcal{C}(H)$.

\paragraph{Stanley's chromatic symmetric function \cite{Sta}:}
For a graph $H$, the \emph{chromatic symmetric function} of $H$ is
$$
X_{{H}}(x)=\sum_{\kappa\in \mathcal{C}({H})}  x_\kappa,
$$ where
\fa{$x_\kappa$ is defined by Equation \ref{eq:xkappa}.}

By definition we have the following.
\begin{proposition} \label{prop:stanleyasGCF}
Let $H$ be a graph. Then
$$X_H(x) = \GX_G(x)$$ \fa{where $G$ is any edge-coloured digraph that has only dashed edges and its underlying graph is isomorphic to $H$.}
\end{proposition} 

\paragraph{Crew-Spirkl's extended chromatic symmetric function \cite{CS}:} A \emph{weighted graph} is a pair  $(H,wt)$ where $H$ is a graph and $wt: V(H)\rightarrow \mathbb{P}$ is a vertex-weight function. The \emph{extended chromatic symmetric function} of $(H,wt)$ is 
$$X_{(H,wt)}(x)=\sum_{\kappa \in \mathcal{C}(H) } x_{\kappa}^{wt},$$ where
$$x_\kappa^{wt}=\prod_{a\in V(H)}x_{\kappa(a)}^{wt(a)}.$$

By definition we have the following.
\begin{proposition}\label{prop:CSasGCF} Let  $(H,wt)$ be a weighted graph. \fathree{For each vertex $a$ in  $V(H)$, consider a cycle $C_a$ isomorphic to $C_{wt(a)}$, and let $v(a)$ be an arbitrary vertex of $C_a$.}
\fa{Define $G$ to be an edge-coloured digraph with the vertex set 
$$\biguplus_{a\in V(H)} V(C_{a})$$ and the edge set
  $$\left(\biguplus_{a\in V(H)} E(C_{a})\right)\uplus \{v(a)\dashrightarrow v(b): ab\in E(H)\}.$$
  Then $$X_{(H,wt)}(x)=\GX_G(x).$$}
  \end{proposition} 

\paragraph{Shareshian-Wachs' chromatic quasisymmetric function \cite{SW}:}
For a graph $H$ with $V(H)$ a subset of $\mathbb{P}$, the \emph{chromatic quasisymmetric function} of $H$ is
$$
X_{{H}}(x,t)=\sum_{\kappa\in \mathcal{C}({H})}t^{{\rm asc}(\kappa)} x_\kappa,
$$
where \fa{$x_\kappa$ is defined by Equation \ref{eq:xkappa} and}  $${\rm asc}(\kappa)=|\{ab\in E(H): a<b \text{~and~} \kappa(a)<\kappa(b)\}|.$$

By definition we have the following.
\begin{proposition}\label{prop:SWasGCF} 
Let $H$ be a graph such that $V(H)$ is a subset of $\mathbb{P}$. Then
$$X_H(x,t)=\GX_G(x,t)$$ where $G$ is any edge-coloured digraph that has only dashed edges with $E(G)=\{(a,b): ab\in E(H), a<b\}$ and its underlying graph is isomorphic to $H$. 
\end{proposition}

\paragraph{Ellzey's chromatic quasisymmetric function \cite{E17}:}
Let $G$ be a digraph.  A \emph{proper} vertex-colouring of $G$ is a vertex-colouring of $G$ such that the colours of adjacent vertices are different. 
Then the \emph{chromatic quasisymmetric function} of $G$ is
$$
Z_{{G}}(x,t)=\sum_{\kappa\in \mathcal{C}({G})}t^{{\rm asc}(\kappa)} x_\kappa,
$$  where \fa{$x_\kappa$ is defined by Equation \ref{eq:xkappa}} \fathree{and $$ {\rm asc}(\kappa)=|\{(a,b)\in E(G):  \kappa(a)<\kappa(b)\}|.$$}

By definition we have the following.
\begin{proposition} \label{prop:EasGCF}
Let $G$ be a digraph. Then
$$Z_{G}(x,t)=\GX_{G}(x,t)$$where $G$ has been edge-coloured with only dashed edges.
\end{proposition}

\begin{remark}
 It is known that Shareshian-Wachs' chromatic quasisymmetric functions of graphs and Ellzey's chromatic quasisymmetric functions of digraphs are symmetric functions when $t=1$. However, the generalized chromatic function $\GX_{G}(x)$ is a quasisymmetric function (not necessarily a symmetric function) as seen in Proposition~\ref{prop:gcfisqsym}. For instance, in Section \ref{bases}, we show that many bases of $\QSym(x)$ can be realized by the generalized chromatic functions of edge-coloured digraphs. 
\end{remark}

Recall that an \emph{orientation} of an undirected graph $H$ is a digraph $G$ with the same vertices,
so that for every edge $ab$ of $H$, exactly one of $(a, b)$ and $(b, a)$ is an edge of $G$.  A \emph{weak cycle} of an orientation $G$ of $H$ is a subgraph $C$ of $G$ such that the underlying graph of $C$, $\overline{C}$, is a cycle in $H$. 
For $k \geq 1$, $G$ is \emph{$k$-balanced}
if for any weak cycle $C$ of $G$ with $E(\overline{C})=\{a_ia_{i+1}: i=1,2,\dots, n-1\}\faa{\uplus} \{ a_1a_{n} \}$, there are at least $k$ edges in $G$ of the form $(a_i,a_{i+1})$ and at least $k$ edges of the form $(a_{i+1},a_i)$ (subscripts are taken modulo $n$).

 Let $\kappa : V (H) \rightarrow \mathbb{P}$ be a proper
vertex-colouring of $H$. Then the orientation induced by $\kappa$ is the orientation $G_\kappa$ where each edge
is directed towards the vertex with the greater colour. If $G_\kappa$ is $k$-balanced, then $\kappa$ is called
a \emph{$k$-balanced colouring}.

\paragraph{Humpert's  $k$-balanced chromatic quasisymmetric function \cite{H}:}
Given a  graph $H$ with $n$ vertices and any positive integer $k$, the \emph{$k$-balanced chromatic quasisymmetric function} of $H$ is 
$$\svwt{X_{H}^k(x)}=\sum x_\kappa,$$ where the sum runs over all $k$-balanced colourings $\kappa: V(H)\rightarrow \mathbb{P}$ and \fa{ $x_\kappa$ is defined by Equation \ref{eq:xkappa}.}

By definition we have the following.
\begin{proposition} \label{prop:HumpasGCF} Let $H$ be a graph and $k$ be a positive integer. Then
$$\svwt{X_{H}^k(x)}=\sum \GX_{G}(x)$$
where the sum runs over all digraphs $G$ such that $G$ is a $k$-balanced orientation of $H$ that has been edge-coloured with only solid edges. 
\end{proposition}

\section{Product and Coproduct Formulas for Generalized Chromatic Functions}\label{prod-coprod}

We now establish the product and coproduct formulas for generalized chromatic functions. The following is straightforward because any pair of independent proper colourings of edge-coloured digraphs $G_1$ and $G_2$  naturally corresponds to a proper colouring of $G_1\uplus G_2$. 
 
\begin{proposition}\label{prod} 
Let $G_1$  and $G_2$ be edge-coloured digraphs.    Then 
\[
\GX_{G_1 }(x,t)\GX_{G_2}(x,t)=\GX_{ {G_1\uplus G_2}}(x,t).
\]
\end{proposition}

To give our coproduct formula, we first need to introduce some notation. Let $G$ be an edge-coloured digraph. An \emph{induced subdigraph} $F=(V(F),E(F))$ of $G$ is a subdigraph of $G$ whose edge set $E(F)$ is the set of all edges $(a,b)$ of $G$ with $a,b\in V(F)$.  
For a subset $A$ of $V(G)$, let $G|_A$ be the induced subdigraph of $G$ with vertex set $A$ where the colours of edges of $G|_A$ are the same as their colours in $G$. 
   An edge-coloured subdigraph $F$  of $G$ is called a \emph{$\{\rightarrow,\Rightarrow\}$-induced subdigraph} of $G$ 
  if $F$ is an induced subdigraph of $G$ \fa{such that} if $a\in V(F)$ and either $a\rightarrow b$ or $a\Rightarrow b$ in $G$, then $b\in V(F)$.

\begin{example}
\fathree{The
$\{\rightarrow,\Rightarrow\}\text{-}$induced subdigraphs of the following edge-coloured digraph
}
$$
  
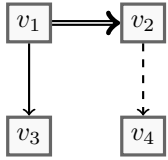
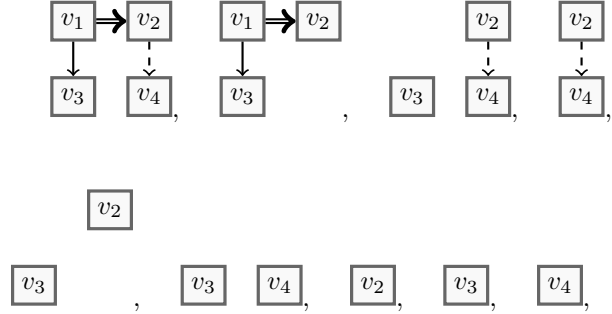
\begin{tikzpicture}[
roundnode/.style={circle, draw=black!60, fill=gray!10, very thick, minimum size=7mm},
squarednode/.style={rectangle, draw=black!60, fill=gray!5, very thick, minimum size=5mm},]
 \node[squarednode](1) at (0,0){$v_1$};
  \node[squarednode](2) at (1.5,0){$v_2$};
    \node[squarednode](3) at (0,-1.5){$v_3$};
       \node[squarednode](4) at (1.5,-1.5){$v_4$};
    \draw[thick,double,->](1)--(2);
       \draw[thick,dashed,->](2)--(4);
        \draw[thick,->](1)--(3);
 \end{tikzpicture}
$$
\fathree{are} 
$$ \begin{tikzpicture}[
roundnode/.style={circle, draw=black!60, fill=gray!10, very thick, minimum size=7mm},
squarednode/.style={rectangle, draw=black!60, fill=gray!5, very thick, minimum size=5mm},]
 \node[squarednode](1) at (0,0){$v_1$};
  \node[squarednode](2) at (1,0){$v_2$};
    \node[squarednode](3) at (0,-1){$v_3$};
       \node[squarednode](4) at (1,-1){$v_4$};
    \draw[thick,double,->](1)--(2);
       \draw[thick,dashed,->](2)--(4);
        \draw[thick,->](1)--(3);
 \end{tikzpicture},~~~~ \begin{tikzpicture}[
roundnode/.style={circle, draw=black!60, fill=gray!10, very thick, minimum size=7mm},
squarednode/.style={rectangle, draw=black!60, fill=gray!5, very thick, minimum size=5mm},]
 \node[squarednode](1) at (0,0){$v_1$};
  \node[squarednode](2) at (1,0){$v_2$};
    \node[squarednode](3) at (0,-1){$v_3$};
    \draw[thick,double,->](1)--(2);
        \draw[thick,->](1)--(3);
 \end{tikzpicture},~~~~
  \begin{tikzpicture}[
roundnode/.style={circle, draw=black!60, fill=gray!10, very thick, minimum size=7mm},
squarednode/.style={rectangle, draw=black!60, fill=gray!5, very thick, minimum size=5mm},]
  \node[squarednode](2) at (1,0){$v_2$};
    \node[squarednode](3) at (0,-1){$v_3$};
      \node[squarednode](4) at (1,-1){$v_4$};
       \draw[thick,dashed,->](2)--(4);
 \end{tikzpicture},~~~~
  \begin{tikzpicture}[
roundnode/.style={circle, draw=black!60, fill=gray!10, very thick, minimum size=7mm},
squarednode/.style={rectangle, draw=black!60, fill=gray!5, very thick, minimum size=5mm},]
  \node[squarednode](2) at (1,0){$v_2$};
      \node[squarednode](4) at (1,-1){$v_4$};
       \draw[thick,dashed,->](2)--(4);
 \end{tikzpicture},$$

  \vspace{3mm}
 
 $$
  \begin{tikzpicture}[
roundnode/.style={circle, draw=black!60, fill=gray!10, very thick, minimum size=7mm},
squarednode/.style={rectangle, draw=black!60, fill=gray!5, very thick, minimum size=5mm},]
  \node[squarednode](2) at (1,0){$v_2$};
    \node[squarednode](3) at (0,-1){$v_3$};
 \end{tikzpicture},~~~~
  \begin{tikzpicture}[
roundnode/.style={circle, draw=black!60, fill=gray!10, very thick, minimum size=7mm},
squarednode/.style={rectangle, draw=black!60, fill=gray!5, very thick, minimum size=5mm},]
    \node[squarednode](3) at (0,-1){$v_3$};
      \node[squarednode](4) at (1,-1){$v_4$};
 \end{tikzpicture},~~~~
  \begin{tikzpicture}[
roundnode/.style={circle, draw=black!60, fill=gray!10, very thick, minimum size=7mm},
squarednode/.style={rectangle, draw=black!60, fill=gray!5, very thick, minimum size=5mm},]
  \node[squarednode](2) at (1,0){$v_2$};
 \end{tikzpicture},~~~~
  \begin{tikzpicture}[
roundnode/.style={circle, draw=black!60, fill=gray!10, very thick, minimum size=7mm},
squarednode/.style={rectangle, draw=black!60, fill=gray!5, very thick, minimum size=5mm},]
    \node[squarednode](3) at (0,-1){$v_3$};
 \end{tikzpicture},  ~~~~
  \begin{tikzpicture}[
roundnode/.style={circle, draw=black!60, fill=gray!10, very thick, minimum size=7mm},
squarednode/.style={rectangle, draw=black!60, fill=gray!5, very thick, minimum size=5mm},]
      \node[squarednode](4) at (1,-1){$v_4$};
 \end{tikzpicture},~~~~
 \emptyset.
$$

\end{example}

We now give the coproduct formula for $\GX_{G}(x)$. 

\begin{theorem}\label{coprod}
Let $G$ be an edge-coloured digraph.  Then 
$$\Delta(\GX_{G}(x))=\sum  \GX_{G|_{V(G)-V(F)}}(x)\otimes \GX_{F}(x),$$ where the sum runs over all $\{\rightarrow,\Rightarrow\}$-induced subdigraphs $F$ of $G$.
\end{theorem} 

\begin{proof}
 Recall that the coproduct of $\QSym(x)$ can be seen as the composite of the following functions
$$
\QSym(x) \cong \QSym(x,y) \rightarrow \QSym(x)\otimes \QSym(y) \cong \QSym(x)\otimes \QSym(x),
$$
where $\QSym(x) \cong \QSym(x,y)$ is defined by $f(x_1,x_2,\ldots)\mapsto f(x_1,x_2,\ldots,y_1,y_2,\ldots)$.
We have $$\GX_{G}(x,y)=\sum_{\kappa} \prod_{a\in V(G)} x_{\kappa(a)},$$ where the sum runs over all functions $\kappa$ from $V(G)$ to the alphabet $$\{x_1<x_2<\cdots <y_1<y_2<\cdots\}$$ such that 
 when $(a,b)$ is in $E(G)$:
\begin{itemize}
\item[(i)] If $a\dashrightarrow b$, then $\kappa(a)\neq \kappa(b)$.
\item[(ii)] If  $a\rightarrow b$, then $\kappa(a) < \kappa(b)$.
\item[(iii)] If  $a\Rightarrow b$, then $\kappa(a)\leq \kappa(b)$.
\end{itemize} 
For the induced subdigraph $F$ with the vertex set $\{a :\kappa(a)\in \{y_1,y_2,\ldots \}\}$, we see that if $a\in V(F)$ and either  $a\rightarrow b$ or $a\Rightarrow b$ in $G$, then $\kappa(a)\leq \kappa(b)$, and so $b\in V(F)$. Therefore, $F$ is a $\{\rightarrow,\Rightarrow\}$-induced subdigraph of $G$.  Also, note that the rest of vertices produce the edge-coloured digraph $G|_{V(G)-V(F)}$. Applying the above composite, we have 
$$\GX_{G}(x)\mapsto \GX_{G}(x,y) \mapsto \sum  \GX_{G|_{V(G)-V(F)}}(x)\otimes \GX_{F}(y)\mapsto \sum \GX_{G|_{V(G)-V(F)}}(x)\otimes \GX_{F}(x),$$ where the  sums run over all $\{\rightarrow,\Rightarrow\}$-induced subdigraphs $F$ of $G$.
\end{proof} 

\begin{remark}
Theorem \ref{coprod} is a generalization of \cite[Proposition 3.4]{AWvW}. Moreover, when $G$ is an edge-coloured digraph that has only dashed edges, then the generalized chromatic function $\GX_{G}(x)$ is equal to Stanley's chromatic symmetric function of the underlying graph of $G$, so this coproduct formula will give a coproduct formula for Stanley's chromatic symmetric functions. 
\end{remark}

\begin{example}
$$\Delta( \GX_{ \begin{tikzpicture}[
roundnode/.style={circle, draw=black!60, fill=gray!10, very thick, minimum size=7mm},
squarednode/.style={rectangle, draw=black!60, fill=gray!5, very thick, minimum size=5mm},]
 \node[squarednode](1) at (0,0){$v_1$};
  \node[squarednode](2) at (1,0){$v_2$};
    \node[squarednode](3) at (0,-1){$v_3$};
       \node[squarednode](4) at (1,-1){$v_4$};
    \draw[thick,double,->](1)--(2);
       \draw[thick,dashed,->](2)--(4);
        \draw[thick,->](1)--(3);
 \end{tikzpicture} }(x))=1\otimes \GX_{ \begin{tikzpicture}[
roundnode/.style={circle, draw=black!60, fill=gray!10, very thick, minimum size=7mm},
squarednode/.style={rectangle, draw=black!60, fill=gray!5, very thick, minimum size=5mm},]
 \node[squarednode](1) at (0,0){$v_1$};
  \node[squarednode](2) at (1,0){$v_2$};
    \node[squarednode](3) at (0,-1){$v_3$};
       \node[squarednode](4) at (1,-1){$v_4$};
    \draw[thick,double ,->](1)--(2);
       \draw[thick,dashed,->](2)--(4);
        \draw[ thick,->](1)--(3);
 \end{tikzpicture} }(x) +
 \GX_{
  \begin{tikzpicture}[
roundnode/.style={circle, draw=black!60, fill=gray!10, very thick, minimum size=7mm},
squarednode/.style={rectangle, draw=black!60, fill=gray!5, very thick, minimum size=5mm},]
      \node[squarednode](4) at (1,-1){$v_4$};
 \end{tikzpicture}
 }(x)\otimes 
 \GX_{\begin{tikzpicture}[
roundnode/.style={circle, draw=black!60, fill=gray!10, very thick, minimum size=7mm},
squarednode/.style={rectangle, draw=black!60, fill=gray!5, very thick, minimum size=5mm},]
 \node[squarednode](1) at (0,0){$v_1$};
  \node[squarednode](2) at (1,0){$v_2$};
    \node[squarednode](3) at (0,-1){$v_3$};
    \draw[thick,double,->](1)--(2);
        \draw[thick,->](1)--(3);
 \end{tikzpicture}}(x) + 
  $$
  
   \vspace{3mm}
  
  $$
  \GX_{
   \begin{tikzpicture}[
roundnode/.style={circle, draw=black!60, fill=gray!10, very thick, minimum size=7mm},
squarednode/.style={rectangle, draw=black!60, fill=gray!5, very thick, minimum size=5mm},]
 \node[squarednode](1) at (0,0){$v_1$};
 \end{tikzpicture}
  }(x)
  \otimes 
   \GX_{\begin{tikzpicture}[
roundnode/.style={circle, draw=black!60, fill=gray!10, very thick, minimum size=7mm},
squarednode/.style={rectangle, draw=black!60, fill=gray!5, very thick, minimum size=5mm},]
  \node[squarednode](2) at (1,0){$v_2$};
    \node[squarednode](3) at (0,-1){$v_3$};
      \node[squarednode](4) at (1,-1){$v_4$};
       \draw[thick,dashed,->](2)--(4);
 \end{tikzpicture}}(x) +
 \GX_{ \begin{tikzpicture}[
roundnode/.style={circle, draw=black!60, fill=gray!10, very thick, minimum size=7mm},
squarednode/.style={rectangle, draw=black!60, fill=gray!5, very thick, minimum size=5mm},]
 \node[squarednode](1) at (0,0){$v_1$};
    \node[squarednode](3) at (0,-1){$v_3$};
        \draw[thick,->](1)--(3);
 \end{tikzpicture}}(x)
 \otimes \GX_{ \begin{tikzpicture}[
roundnode/.style={circle, draw=black!60, fill=gray!10, very thick, minimum size=7mm},
squarednode/.style={rectangle, draw=black!60, fill=gray!5, very thick, minimum size=5mm},]
  \node[squarednode](2) at (1,0){$v_2$};
      \node[squarednode](4) at (1,-1){$v_4$};
       \draw[thick,dashed,->](2)--(4);
 \end{tikzpicture}}(x)+
 \GX_{  \begin{tikzpicture}[
roundnode/.style={circle, draw=black!60, fill=gray!10, very thick, minimum size=7mm},
squarednode/.style={rectangle, draw=black!60, fill=gray!5, very thick, minimum size=5mm},]
 \node[squarednode](1) at (0,0){$v_1$};
      \node[squarednode](4) at (1,-1){$v_4$};
 \end{tikzpicture}}(x)
\otimes \GX_{  \begin{tikzpicture}[
roundnode/.style={circle, draw=black!60, fill=gray!10, very thick, minimum size=7mm},
squarednode/.style={rectangle, draw=black!60, fill=gray!5, very thick, minimum size=5mm},]
  \node[squarednode](2) at (1,0){$v_2$};
    \node[squarednode](3) at (0,-1){$v_3$};
 \end{tikzpicture}} (x)+
 $$
 
 \vspace{3mm}
 
 $$
 \GX_{
  \begin{tikzpicture}[
roundnode/.style={circle, draw=black!60, fill=gray!10, very thick, minimum size=7mm},
squarednode/.style={rectangle, draw=black!60, fill=gray!5, very thick, minimum size=5mm},]
 \node[squarednode](1) at (0,0){$v_1$};
  \node[squarednode](2) at (1,0){$v_2$};
    \draw[thick,double,->](1)--(2);
 \end{tikzpicture}
 }(x) \otimes 
 \GX_{ \begin{tikzpicture}[
roundnode/.style={circle, draw=black!60, fill=gray!10, very thick, minimum size=7mm},
squarednode/.style={rectangle, draw=black!60, fill=gray!5, very thick, minimum size=5mm},]
    \node[squarednode](3) at (0,-1){$v_3$};
      \node[squarednode](4) at (1,-1){$v_4$};
 \end{tikzpicture}}(x)+
 \GX_{
  \begin{tikzpicture}[
roundnode/.style={circle, draw=black!60, fill=gray!10, very thick, minimum size=7mm},
squarednode/.style={rectangle, draw=black!60, fill=gray!5, very thick, minimum size=5mm},]
 \node[squarednode](1) at (0,0){$v_1$};
    \node[squarednode](3) at (0,-1){$v_3$};
      \node[squarednode](4) at (1,-1){$v_4$};
        \draw[thick,->](1)--(3);
 \end{tikzpicture}
 }(x)
 \otimes 
 \GX_{ \begin{tikzpicture}[
roundnode/.style={circle, draw=black!60, fill=gray!10, very thick, minimum size=7mm},
squarednode/.style={rectangle, draw=black!60, fill=gray!5, very thick, minimum size=5mm},]
  \node[squarednode](2) at (1,0){$v_2$};
 \end{tikzpicture}}(x)+
 \GX_{
  \begin{tikzpicture}[
roundnode/.style={circle, draw=black!60, fill=gray!10, very thick, minimum size=7mm},
squarednode/.style={rectangle, draw=black!60, fill=gray!5, very thick, minimum size=5mm},]
 \node[squarednode](1) at (0,0){$v_1$};
  \node[squarednode](2) at (1,0){$v_2$};
         \node[squarednode](4) at (1,-1){$v_4$};
    \draw[thick,double,->](1)--(2);
      \draw[ thick,dashed,->](2)--(4);
 \end{tikzpicture}
 }\otimes 
 \GX_{ \begin{tikzpicture}[
roundnode/.style={circle, draw=black!60, fill=gray!10, very thick, minimum size=7mm},
squarednode/.style={rectangle, draw=black!60, fill=gray!5, very thick, minimum size=5mm},]
    \node[squarednode](3) at (0,-1){$v_3$};
 \end{tikzpicture}}(x)+
  $$
  
  \vspace{3mm}
  
  $$
  \GX_{ \begin{tikzpicture}[
roundnode/.style={circle, draw=black!60, fill=gray!10, very thick, minimum size=7mm},
squarednode/.style={rectangle, draw=black!60, fill=gray!5, very thick, minimum size=5mm},]
 \node[squarednode](1) at (0,0){$v_1$};
  \node[squarednode](2) at (1,0){$v_2$};
    \node[squarednode](3) at (0,-1){$v_3$};
    \draw[thick,double,->](1)--(2);
      \draw[thick,->](1)--(3);
 \end{tikzpicture}}(x)\otimes 
 \GX_{ \begin{tikzpicture}[
roundnode/.style={circle, draw=black!60, fill=gray!10, very thick, minimum size=7mm},
squarednode/.style={rectangle, draw=black!60, fill=gray!5, very thick, minimum size=5mm},]
      \node[squarednode](4) at (1,-1){$v_4$};
 \end{tikzpicture}}(x)+ 
 \GX_{ \begin{tikzpicture}[
roundnode/.style={circle, draw=black!60, fill=gray!10, very thick, minimum size=7mm},
squarednode/.style={rectangle, draw=black!60, fill=gray!5, very thick, minimum size=5mm},]
 \node[squarednode](1) at (0,0){$v_1$};
  \node[squarednode](2) at (1,0){$v_2$};
    \node[squarednode](3) at (0,-1){$v_3$};
       \node[squarednode](4) at (1,-1){$v_4$};
    \draw[thick,double,->](1)--(2);
       \draw[thick,dashed,->](2)--(4);
        \draw[thick,->](1)--(3);
 \end{tikzpicture} }(x)\otimes 1.
  $$

\end{example}

\section{Bases for $\QSym^r(x)$  Using Generalized Chromatic Functions}\label{bases}

In this section, we realize different bases for the Hopf algebras in the following chain as special cases of the generalized chromatic functions of certain edge-coloured digraphs.
$${\rm QSym}(x)={\rm QSym}^1(x)  \supset {\rm QSym}^2(x)  \supset \cdots \supset {\rm QSym}^\infty(x)={\rm Sym}(x)$$

Given a partition $\lambda=(\lambda_1,\lambda_2,\ldots,\lambda_l)\vdash n$,  the edge-coloured digraph $G_{\lambda}$ is an edge-coloured digraph such that:
\begin{itemize}
\item[(i)] The vertex set of ${G}$ is 
$$\{ (i,j): 1\leq i\leq l, 1\leq j \leq \lambda_i\}.$$
\item[(ii)]  $(i,j)\rightarrow (i',j')$ in $G$ if and only if \fa{$i'=i+1$ and $j'=j$.}
\item[(iii)] $(i,j)\Rightarrow (i',j')$ in $G$ if and only if \fa{$i'=i$ and  $j'=j+1$.}
 \end{itemize}

\sxl{Let $\alpha=(\alpha_1,\alpha_2,\dots,\alpha_{k})\vDash[n]$ and $\lambda=(\lambda_1,\lambda_2,\dots, \lambda_{l})=(n^{m_n}, \fss{(n-1)}^{m_{n-1}}, \ldots, 1^{m_1})\vdash n$.} In the following table, we see that many well-known bases of $\Sym(x)$ and $\QSym(x)$ are the generalized chromatic functions of some edge-coloured digraphs produced by the actions of the operators 
$$\begin{tikzpicture} 
\node at (0,0){$\biguplus,$}; 
\node at (1,0){ $\arx$};
\node at (1.3,-0.15){,};
\node at (2,0){ $\ary$};
\node at (2.3,-0.15){,};
\node at (3,0){ $\arz$};
\end{tikzpicture} $$
on the edge-coloured digraphs in Table \ref{ecd}. For these well-known bases, the result follows immediately by definition, and hence readers unfamiliar with the classical definitions may take these results to be the definitions, or refer to \cite{Sta1}. For the upper-fundamental basis of $\QSym(x)$ take the commutative image of the noncommutative upper-fundamental basis in \cite{UpFun}.

\begin{table}[H]
\caption{Bases for $\Sym(x)$ and $\QSym(x)$ reinterpreted}
\label{table:basesQS}
\centering
\begin{tabsize}
\begin{tabular}{ |c|c|c| }
\toprule
\text{{Basis}} &$G$& { $\GX_{G}$ }\\
\midrule
  \begin{tikzpicture} \node at (0,0){Monomial basis of $\Sym(x)$};\end{tikzpicture}   &$ \begin{tikzpicture} 
  \node at (-1.2,0){$\arx_{j=1}^n$};
\node at (-0.6,0){${\Big (}$};
\node at (0,0){$\ary_{i=1}^{m_j}$};
\node at (0.5,0){$~~~~C_{j}$};
\node at (1.1,0){$\Big )$};
\end{tikzpicture}$& \begin{tikzpicture}\node at (0,0){$ m_\lambda$};  \end{tikzpicture} \\
\midrule
  \begin{tikzpicture} \node at (0,0){Augmented monomial basis of $\Sym(x)$};\end{tikzpicture}     & $\begin{tikzpicture} 
\node at (0,0){$\arx_{i=1}^l$};
\node at (0.5,0){$~~~~C_{\lambda_i}$};
\end{tikzpicture}$&\begin{tikzpicture}\node at (0,0){\fathree{$\widetilde{m}_\lambda=\lambda^! m_\lambda$}};  \end{tikzpicture}  \\
\midrule
 \text{Elementary basis of}~$\Sym(x)$  & $\biguplus_{i=1}^l P_{\lambda_i}$& $e_\lambda$\\
\midrule
\text{Augmented elementary basis of}~$\Sym(x)$ & $\biguplus_{i=1}^l K_{\lambda_i}$& $\lambda! e_\lambda$\\
\midrule
 \text{Complete homogeneous basis of}~$\Sym(x)$ & $\biguplus_{i=1}^l Q_{\lambda_i}$& $h_\lambda$ \\
\midrule
 \text{Power sum basis of}~$\Sym(x)$  & $\biguplus_{i=1}^l C_{\lambda_i}$& $p_\lambda$\\
\midrule
 \text{Schur basis of}~$\Sym(x)$  & $G_\lambda$& $s_\lambda$\\
\midrule
 \begin{tikzpicture} \node at (0,0){Monomial basis of $\QSym(x)$};\end{tikzpicture}   & $\begin{tikzpicture} 
\node at (0,0){$\ary_{i=1}^k$};
\node at (0.5,0){$~~~~C_{\alpha_i}$};
\end{tikzpicture}$&\begin{tikzpicture}\node at (0,0){$M_\alpha$};  \end{tikzpicture}\\
\midrule
 \begin{tikzpicture} \node at (0,0){Fundamental basis of $\QSym(x)$};\end{tikzpicture}  & $\begin{tikzpicture} 
\node at (0,0){$\ary_{i=1}^k$};
\node at (0.6,0){$~~~Q_{\alpha_i}$};
\end{tikzpicture}$& \begin{tikzpicture}\node at (0,0){$F_\alpha$};  \end{tikzpicture}  \\
\midrule
 \begin{tikzpicture} \node at (0,0){Upper-fundamental basis of $\QSym(x)$};\end{tikzpicture}    & $\begin{tikzpicture} 
\node at (0,0){$\arz_{i=1}^k$};
\node at (0.5,0){$~~~~C_{\alpha_i}$};
\end{tikzpicture}$& \begin{tikzpicture}\node at (0,0){$\overline{F}_\alpha$};  \end{tikzpicture}  \\
\bottomrule
\end{tabular}
\end{tabsize}
\end{table}

\begin{remark}
If we generalize the edge-coloured digraph $G_\lambda$ for a partition $\lambda$ to $G_\alpha$ for a composition $\alpha$ in the natural way, but restrict the second condition to only the first column, then by definition we have
$$
\GX_{G_\alpha}(x) = \mathfrak{S} ^\ast _{\alpha},
$$ where $\mathfrak{S} ^\ast _{\alpha}$ is the dual immaculate function indexed by $\alpha$ \cite{BBSSZ}. Switching the $\rightarrow$ for $\Rightarrow$ and vice versa in $G_\alpha$ we obtain the row-strict dual immaculate function $\mathcal{R}\mathfrak{S} ^\ast _{\alpha}$ \cite{NSVWVW}.
\end{remark}

Let $G$ be an edge-coloured digraph. Recall that $P(G,p)$  is the number of ways that one can properly colour the vertices of  $G$ with $p$ colours. The following theorem shows that $P(G,p)$ is a polynomial in $p$.

\begin{theorem}\label{thm:gcp}
For any edge-coloured digraph $G$, $P(G,p)$ is a polynomial in $p$. 
\end{theorem} 

\begin{proof}
\fa{Note that when we set the first $p$ variables of $\GX_G(x)$ equal to $1$ and all others to $0$, we have $P(G,p)$. Since $\GX_G(x)$ is a quasisymmetric function by Proposition~\ref{prop:gcfisqsym}, it is a linear combination of the monomial basis elements of $\QSym(x)$. Thus, we only need to show that  for any composition $\alpha=(\alpha_1,\alpha_2,\dots,\alpha_k)$,  
$$M_\alpha(\underbrace{1,1,\dots, 1}_{\text{$p$ times}}, 0,0, \dots)$$
is a polynomial in $p$. Note that this is equal to the number of monomials appearing in $M_{\alpha}(x_1,x_2,\dots,x_p)$. The number of these monomials is the number of ways of choosing $k$ of the $p$ variables, which is ${p \choose k}$, a polynomial in $p$.}
\end{proof}

 In \cite{CvW}, Cho and van Willigenburg constructed an infinite family of bases for $\Sym(x)$ using chromatic symmetric functions of graphs. In the following theorem, we construct an infinite family of bases for $\QSym(x)$ using generalized chromatic functions. 

\begin{theorem}
 Let $\mathcal{F}_i$ be an edge-coloured digraph with $i$ vertices that has only double edges. For each $\alpha=(\alpha_1,\alpha_2,\dots, \alpha_{k})\vDash n$, define
$$
\begin{tikzpicture} 
\node (a) at (0,0){$\mathcal{F}_{\alpha}=$};
\node (b) at (1,0){$\mathcal{F}_{\alpha_1}$};
\node (c) at (1.7,0){$\ary$};
\node (d) at (2.5,0){$\mathcal{F}_{\alpha_2}$};
\node (e) at (3.2,0){$\ary$};
\node (f) at (4,0){$\dots$};
\node (g) at (4.7,0){$\ary$};
\node (h) at (5.6,0){$\mathcal{F}_{\alpha_{k}}.$};
\end{tikzpicture} 
$$ Then  
$$\{\GX_{\mathcal{F}_\alpha} (x):\alpha\vDash n\}$$ is a basis for $\QSym_n(x)$.  
\end{theorem}

\begin{proof}
Writing $\GX_{\mathcal{F}_\alpha}(x)$ in terms of the monomial basis elements of $\QSym(x)$, we have $$\GX_{\mathcal{F}_\alpha}(x)=M_\alpha+\sum_{\beta >\alpha}c_{\alpha,\beta}M_\beta$$ for some coefficients $c_{\alpha,\beta}$. Since each $\GX_{\mathcal{F}_\alpha }(x)$ has a unique leading term $M_\alpha$ under the $<$ order, we can conclude that $\{\GX_{\mathcal{F}_\alpha}(x): \alpha\vDash n\}$ is a basis for $\QSym_n(x)$.   
\end{proof}
\sxl{\begin{example}
If $\mathcal{F}_i=C_{i}$, then $\{ \GX_{\mathcal{F}_\alpha}(x)\}$ is the monomial basis of $\QSym(x)$. If $\mathcal{F}_i=Q_{i}$, then $\{ \GX_{\mathcal{F}_\alpha}(x)\}$ is the fundamental basis of $\QSym(x)$. 
\end{example}}

\fa{We now define a product on $\mathbb{Q}[[x_1,x_2,\dots ]]$  that plays a crucial role in the rest of this section.  Given a multiset $I$ of positive integers, define 
$$x_I=\prod_{i\in I} x_i.$$ Define a product on $\mathbb{Q}[[x_1,x_2,\dots ]]$ by bilinearly extending the operator
$$x_I\odot x_J=\begin{cases} 
x_Ix_J & \text{if~}I\cap J=\emptyset,\\
0 & \text{otherwise.} 
\end{cases} $$
Then $$\GX_{
\begin{tikzpicture} 
\node (b) at (0,0){$G_1$};
\node (c) at (0.5,0){$\arx$};
\node (d) at (1.1,0){$G_2$};
\end{tikzpicture}}(x)=\GX_{G_1}(x)\odot \GX_{G_2}(x).$$}

Let $r$ be a positive integer or infinity.  We now establish several new bases for  $\mathrm{QSym}^r(x)$. 
Let $(\beta,\mu)=((\beta_1,\beta_2,\ldots,\beta_k),(\mu_1,\mu_2,\ldots,\mu_l))$ be an $r$-composition. \fa{By definition we have that
$$
M_{(\beta,\mu)}=\GX_{
\begin{tikzpicture} 
\node at (-0.6,0){${\Big (}$};
\node at (0,0){$\ary_{i=1}^k$};
\node at (0.5,0){$~~~~C_{\beta_i}$};
\node at (1.1,0){$\Big )$};
\node at (1.5,0){$\arx$};
\node at (1.9,0){$\Big ($};
\node at (2.5,0){$\arx_{j=1}^l$};
\node at (3,0){$~~~~C_{\mu_j}$};
\node at (3.6,0){$\Big )$};
\end{tikzpicture}
}(x)=M_\beta\odot \widetilde{m}_\mu.
$$
Define $$S_{(\beta,\mu)}= \GX_{
\begin{tikzpicture} 
\node at (-0.6,0){${\Big (}$};
\node at (0,0){$\ary_{i=1}^k$};
\node at (0.5,0){$~~~~C_{\beta_i}$};
\node at (1.1,0){$\Big )$};
\node at (1.5,0){$\arx$};
\node at (2.1,0){$G_\mu$};
\end{tikzpicture}
}(x)=M_\beta\odot s_\mu,
$$
$$\overline{F}_{(\beta,\mu)}=\GX_{
\begin{tikzpicture} 
\node at (-0.6,0){${\Big (}$};
\node at (0,0){$\arz_{i=1}^k$};
\node at (0.5,0){$~~~~C_{\beta_i}$};
\node at (1.1,0){$\Big )$};
\node at (1.5,0){$\arx$};
\node at (1.9,0){$\Big ($};
\node at (2.5,0){$\arx_{j=1}^l$};
\node at (3,0){$~~~~C_{\mu_j}$};
\node at (3.6,0){$\Big )$};
\end{tikzpicture}}(x)=\overline{F}_{\beta}\odot\widetilde{m}_\mu,$$
and
$$\overline{S}_{(\beta,\mu)}=\GX_{
\begin{tikzpicture} 
\node at (-0.6,0){${\Big (}$};
\node at (0,0){$\arz_{i=1}^k$};
\node at (0.5,0){$~~~~C_{\beta_i}$};
\node at (1.1,0){$\Big )$};
\node at (1.5,0){$\arx$};
\node at (2.1,0){$G_\mu$};
\end{tikzpicture}
}(x)=\overline{F}_{\beta}\odot s_\mu.$$
}
Let $K_{\mu,\nu}$ be the Kostka number indexed by partitions $\mu$ and $\nu$ \cite{Sta1}.
\fa{\begin{theorem}\label{thm:r-bases}
Let $(\beta, \mu)$ be an $r$-composition. Then we have that:
\begin{itemize} 
\item[(i)] $\displaystyle{S_{(\beta,\mu)}=\sum_{\nu\preceq \mu} \frac{K_{\mu,\nu}}{\nu^!}M_{(\beta,\nu)}}$.
\item[(ii)] 
$\displaystyle{ \overline{F}_{(\beta,\mu)}=\sum_{\gamma\leq \beta} M_{(\gamma,\mu)}}$.
\item[(iii)] 
$\displaystyle{ \overline{S}_{(\beta,\mu)}=\sum_{\gamma\leq \beta,  \nu\preceq \mu}  \frac{K_{\mu,\nu}}{\nu^!}M_{(\gamma,\nu)}}.$
\end{itemize} 
\end{theorem} 
\begin{proof}
Note that $$s_\mu=\sum_{\nu\preceq \mu} K_{\mu,\nu} \frac{\widetilde{m}_{\nu}}{\nu^!}\quad \text{and} \quad \overline{F}_\beta=\sum_{\gamma\leq \beta} M_\gamma.$$
Therefore, 
$$S_{(\beta,\mu)}=M_{\beta}\odot s_\mu=\sum_{\nu\preceq \mu} \frac{K_{\mu,\nu}}{\nu^!}\left( M_\beta\odot \widetilde{m}_\nu\right)=\sum_{\nu\preceq \mu} \frac{K_{\mu,\nu}}{\nu^!}M_{(\beta,\nu)},$$
$$ \overline{F}_{(\beta,\mu)}=\overline{F}_\beta \odot \widetilde{m}_\mu=\sum_{\gamma\leq \beta} M_\gamma \odot \widetilde{m}_\mu=\sum_{\gamma\leq \beta} M_{(\gamma,\mu)},$$
and 
$$ \overline{S}_{(\beta,\mu)}= \overline{F}_\beta \odot s_\mu=\sum_{\gamma\leq \beta, \nu\preceq \mu}  \frac{K_{\mu,\nu}}{\nu^!}\left(M_\gamma \odot \widetilde{m}_\nu\right)=\sum_{\gamma\leq \beta,  \nu\preceq \mu}  \frac{K_{\mu,\nu}}{\nu^!} M_{(\gamma,\nu)}.$$
\end{proof} 
}
\fa{Since $K_{\mu,\mu}=1$, for any partition $\mu$, we have the following corollary. 
\begin{corollary}\label{coro:r-bases} Each of the following is a basis for $\QSym_n^r(x)$.
\begin{itemize}
\item[(i)] $\{ {S}_{(\beta,\mu)}: (\beta,\mu) \text{~is an $r$-composition of $n$} \}$.
\item[(ii)] $\{ \overline{F}_{(\beta,\mu)}: (\beta,\mu) \text{~is an $r$-composition of $n$} \}$.
\item[(iii)] $\{ \overline{S}_{(\beta,\mu)}: (\beta,\mu) \text{~is an $r$-composition of $n$} \}$.
\end{itemize} 
\end{corollary}
}

\begin{remark} 
\fathree{The sets 
\[\{ F_\beta \odot \widetilde{m}_\mu: (\beta,\mu) \text{~is an $r$-composition of $n$} \}\]
and 
\[\{ F_\beta \odot s_\mu: (\beta,\mu) \text{~is an $r$-composition of $n$} \}\]
 are not bases for $\QSym^r_n(x)$ since we do not necessarily have $F_\beta \odot \widetilde{m}_\mu, F_\beta \odot s_\mu \in \QSym^r_n(x)$. As an example,  $((2,2),\emptyset)$ is a $2$-composition, and  $$F_{(2,2)}\odot \widetilde{m}_\emptyset=F_{(2,2)}\odot s_\emptyset=F_{(2,2)}=M_{(2,2)}+M_{(2,1,1)}+M_{(1,1,2)}+M_{(1,1,1,1)}\not\in \QSym_4^2(x).$$}
\end{remark}

\section{{Symmetric Functions and Generalizations} in Noncommuting Variables}\label{Hopf algebras noncomm}

This section introduces the Hopf algebras of symmetric, quasisymmetric, and $r$-quasisymmetric functions in noncommuting variables. The bases of these Hopf algebras are indexed by set partitions, set compositions, and $r$-set-compositions, respectively. We recall the definitions and notation related to these combinatorial objects.

\subsection{Set partitions, set compositions, and $r$-set-compositions}

A \emph{set partition} $\Pi$ of a finite set $A$  is a set consisting of mutually disjoint nonempty subsets $\Pi_1,\Pi_2,\ldots,\Pi_l$ of $A$ such that their union is $A$; this is denoted by $\Pi=\Pi_1/\Pi_2/\cdots/\Pi_l\vdash A$. 
 Each $\Pi_i$ is called a \emph{block} of the set partition $\Pi$, and the \emph{length} of $\Pi$ is $l$. By convention, we denote by $\emptyset$ the unique empty set partition of $[0]=\emptyset$.
 Let ${\lambda}(\Pi)=(|\Pi_1|,|\Pi_2|,\ldots,|\Pi_l|)$  where we assume that $|\Pi_1|\geq |\Pi_2|\geq \cdots \geq |\Pi_l|$. We say $\Pi$ is of \emph{shape} $\lambda$ if $\lambda(\Pi)=\lambda$. \fa{The \emph{standardization} of $\Pi$, ${\rm std}(\Pi)$, is the set partition of $[|A|]$ such that the $i$th smallest element of $A$ is replaced by $i$ in each block of $\Pi$ for all $i\in[|A|]$.}
  For example, $\Pi=35/67/9\vdash \{3,5,6,7,9\}$ has length 3, $\lambda(\Pi)=(2,2,1)$, and  ${\rm std}(\Pi)=12/34/5$.

A \emph{set composition} $\Phi$ of a finite set $A$ is a list  of mutually disjoint  nonempty subsets  $\Phi_1,\Phi_2,\ldots,\Phi_k$ of $A$ such that their union is $A$; this is denoted by $(\Phi_1|\Phi_2|\cdots|\Phi_k)\vDash A$.
 Each $\Phi_i$ is called a \emph{block} of the set composition $\Phi$, and the \emph{length} of $\Phi$ is $k$.  By convention, we denote by $\emptyset$ the unique empty set composition of $[0]=\emptyset$. Let $\alpha(\Phi)=(|\Phi_1|,|\Phi_2|,\ldots,|\Phi_k|)$. We say $\Phi$ is of \emph{shape} $\alpha$ if $\alpha(\Phi)=\alpha$. \fa{The \emph{standardization} of $\Phi$, ${\rm std}(\Phi)$, is the set composition of $[|A|]$ such that the $i$th smallest element of $A$ is replaced by $i$ in each block of $\Phi$ for all $i\in[|A|]$.}
  For example, $\Phi=(35|9|67)\vDash\{3,5,6,7,9\}$ has length 3, $\alpha(\Phi)=(2,1,2)$, and ${\rm std}(\Phi)=(12|5|34)$.
\fa{We write the elements within each block in increasing order and consider them ordered.}
\svwt{We say a set composition $\Psi$ \emph{corrupts} $\Phi$  if  $\Psi$ is  $\Phi$ with some bars removed, and say that $\Phi$ \emph{reforms} $\Psi$ if $\Psi$ is $\Phi$ with some bars added. In particular, the numbers of both {\emph{must}} be written in the same order. For example, $(13|2|456)$ corrupts $(13|2|4|56)$, and $(1|3|2|4|56)$ reforms $(13|2|4|56)$  but  $(3|1|2|4|56)$ does not reform $(13|2|4|56)$. }

Let $r$ be a positive integer or infinity. \fa{Let $A\subseteq B \subset \mathbb{P}$ with $B$ finite.} A pair $(\Phi,\Pi)$ where $\Phi$ is a set composition of $A$ and $\Pi$ is a set partition of $B- A$ is called an \emph{$r$-set-composition} of $B$  if $(\alpha(\Phi),\lambda(\Pi))$ is an $r$-composition.  
\fa{The \emph{standardization} of an $r$-set-composition $(\Phi,\Pi)$ of $B$, ${\rm std}(\Phi,\Pi)=(\Psi,\Omega)$, is the $r$-set-composition of $B$ such that  the $i$th smallest element of $B$ is replaced by $i$ in each block of $\Phi$ and $\Pi$ for all $i\in[|B|]$.
}
\sxl{For example, $((36|29), 4/5/8)$ is a $2$-set-composition of $\{2,3,4,5,6,8,9\}$ and we have ${\rm std}((36|29), 4/5/8)=((25|17), 3/4/6)$.}

\subsection{Quasisymmetric functions in noncommuting variables}
The Hopf algebra of quasisymmetric functions in noncommuting variables first appeared in \cite{Hi04}. It is realized as power series in noncommuting variables, and the bases of this Hopf algebra are indexed by set compositions. 

Let $\mathbb{Q}\langle \langle \bx_1,\bx_2,\ldots \rangle \rangle$ be the algebra of formal power series in infinitely many noncommuting variables \fathree{$\bx=\{\bx_1,\bx_2,\ldots\}$} over $\mathbb{Q}$.

\begin{definition}
A \emph{quasisymmetric function in noncommuting variables} is a formal power series $f\in \mathbb{Q}\langle \langle \bx_1,\bx_2,\ldots \rangle \rangle$ such that:
\begin{itemize}
\item[(i)] The degrees of the monomials in $f$ are bounded.
\item[(ii)] For every set composition $\Phi=(\Phi_1|\Phi_2|\cdots|\Phi_k)\vDash [n]$, all monomials $\bx_{i_1}\bx_{i_2} \cdots \bx_{i_n}$ in $f$ satisfying the following conditions have the same coefficient. 
\sxl{\begin{itemize}
\item[(a)] $i_j=i_\ell$ if $j$ and $\ell$ are in the same block of $\Phi$.
\item[(b)] $i_j<i_\ell$ if $j\in \Phi_p$ and $\ell\in \Phi_q$ with $p<q$.
\end{itemize}
}
\end{itemize}

The set of all quasisymmetric functions in noncommuting variables is denoted by $\NCQSym(\bx)$. 
\end{definition}

The vector space $\NCQSym(\bx)$ is  a Hopf algebra where its product is the same as the product of the formal power series in noncommuting variables, and its coproduct $\Delta$ is defined as follows (for more details see \cite[Section 5]{BZ}). Evaluate an element $f(\bx)\in \NCQSym(\bx)\cong \NCQSym(\bx,{\bf y})$ as $f(\bx,{\bf y})$ using the linearly ordered noncommuting variables $(\bx, {\bf y})=(\bx_1<\bx_2< \cdots <{\bf y}_1< {\bf y}_2<\cdots)$. Denote by $\bar{f}(\bx,\by)$  the image of $f(\bx,{\bf y})$ after imposing the partial commutativity relations
$$\bx_i {\bf y}_j=\by_j\bx_i \quad \text{for every pair~} (\bx_i,\by_j)\in \bx\times \by.$$
We have that $\bar{f}(\bx,\by)$ lies in a subalgebra isomorphic to 
$\NCQSym(\bx)\otimes \NCQSym(\by)$. Let the image of $\bar{f}(\bx,\by)$ in $\NCQSym(\bx)\otimes \NCQSym(\by)$ be $$\sum \bar{f}_1(\bx)\otimes \bar{f}_2(\by).$$ 
Applying the following isomorphism
$$\NCQSym(\bx)\otimes \NCQSym(\by) \cong \NCQSym(\bx)\otimes \NCQSym(\bx),$$ we have
$$\sum \bar{f}_1(\bx)\otimes \bar{f}_2(\by)\mapsto \sum \bar{f}_1(\bx)\otimes \bar{f}_2(\bx).$$
 Now we define $\Delta: \NCQSym(\bx) \rightarrow \NCQSym(\bx)\otimes \NCQSym(\bx)$  such that 
$$
\Delta(f(\bx))= \sum \bar{f}_1(\bx)\otimes \bar{f}_2(\bx).
$$

\subsection{Symmetric functions in noncommuting variables}

The Hopf algebra of symmetric functions in noncommuting variables is a Hopf subalgebra of $\NCQSym(\bx)$.
\begin{definition}
A \emph{symmetric function in noncommuting variables} is a formal power series $f\in \mathbb{Q}\langle\langle \bx_1,\bx_2,\ldots \rangle \rangle$  such that:
\begin{itemize} 
\item[(i)] The degrees of the monomials in $f$ are bounded. 
\item[(ii)] For any permutation $\sigma\in \mathfrak{S}_\infty$, $$\sigma \cdot f(\bx_1,\bx_2,\ldots)=f(\bx_{\sigma(1)},\bx_{\sigma(2)}, \ldots)=f(\bx_1,\bx_2,\ldots).$$
\end{itemize} 
 The set of all symmetric functions in noncommuting variables is denoted by $\NCSym(\bx)$. 
 \end{definition}

\subsection{$r$-quasisymmetric functions in noncommuting variables}

 We now introduce the new Hopf algebra of $r$-quasisymmetric functions in noncommuting variables. 

For each $r$-set-composition $(\Phi,\Pi)$ of $[n]$, define the \emph{$r$-dominant monomial function} in noncommuting variables to be 
$${\bf M}_{(\Phi,\Pi)}=\sum_{(i_1,i_2,\ldots,i_n)} \bx_{i_1}\bx_{i_2}\cdots \bx_{i_n}$$ where the sum runs over all tuples $(i_1,i_2,\ldots,i_n)$ such that:
\begin{itemize}
\item[(i)] $i_j=i_\ell$ if and only if $j$ and $\ell$ are  in the same block of either $\Phi$ or $\Pi$.
\item[(ii)] $i_j<i_\ell$ if $j\in \Phi_p$ and $\ell\in \Phi_q$ with $p<q$.
\end{itemize}

Let
$$
{\rm NCQSym}^r(\bx)=\bigoplus_{n\geq 0}{\rm NCQSym}_n^r(\bx),
$$
where
$${\rm NCQSym}_n^r(\bx)=\mathbb{Q}\text{-span}\{ {\bf M}_{(\Phi,\Pi)}: (\Phi,\Pi) \text{~ is an $r$-set-composition of $[n]$}\}.$$
We have that 
$${\rm NCQSym}(\bx)={\rm NCQSym}^1(\bx)  \supset {\rm NCQSym}^2(\bx)  \supset \cdots \supset {\rm NCQSym}^\infty(\bx)={\rm NCSym}(\bx).$$

We will show later in Section~\ref{sec:basesR} that ${\rm NCQSym}^r(\bx)$ is a Hopf algebra and give natural bases for the Hopf algebras in this chain realized as generalized chromatic functions in noncommuting variables, which are defined in the next section.  

\section{Generalized Chromatic Functions in Noncommuting Variables}\label{labelled}
A \emph{labelled edge-coloured digraph} is an edge-coloured digraph whose vertex set is a subset of $\mathbb{P}$.   We usually denote a labelled edge-coloured digraph by $\bG$. 
\begin{figure}[H]
\begin{center}
  \begin{tikzpicture}[
roundnode/.style={circle, draw=black!60, fill=gray!10, very thick, minimum size=7mm},
squarednode/.style={rectangle, draw=black!60, fill=gray!5, very thick, minimum size=5mm},]
 \node[squarednode](1) at (1,1){$1$};
  \node[squarednode](3) at (2,0){$6$};
    \node[squarednode](2) at (3,1){$3$};
    \draw[thick, dashed,->](1)--(2);
       \draw[ thick,->](2)--(3);
        \draw[ thick,double,->](3)--(1);
 \end{tikzpicture}
   \caption{A labelled edge-coloured digraph}
   \label{ex2}
\end{center}
\end{figure}
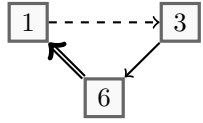

We frequently use the labelled edge-coloured digraphs in the following table.

\begin{table}[H]
\caption{Useful labelled edge-coloured digraphs}
 \label{lecd}
\centering
\begin{tabsize}
\begin{tabular}{|c|c|}
\toprule
 \text{Notation} & \text{Expression}  \\
\midrule
 $C_A$ &   \text{$C_{|A|}$ with vertex set $A$ }\\
\midrule
 $P_A$ & \text{$P_{|A|}$ with vertex set $A$  such that if $a\rightarrow b$, then $a<_{\mathbb{Z}}b$}\\
\midrule
 $Q_A$ &  \text{$Q_{|A|}$ with vertex set  $A$ such that if $a\Rightarrow b$, then $a<_{\mathbb{Z}} b$}\\
\midrule
 $K_A$ & \text{$K_{|A|}$ with vertex set $A$ such that if $a\dasharrow b$, then $a<_{\mathbb{Z}} b$}\\
\bottomrule
\end{tabular}
\end{tabsize}
\end{table}

\begin{figure}[H]
 \begin{center}
   \begin{tikzpicture}[
roundnode/.style={circle, draw=black!60, fill=gray!10, very thick, minimum size=7mm},
squarednode/.style={rectangle, draw=black!60, fill=gray!5, very thick, minimum size=5mm},]
 \node[squarednode](1) at (1.3,0){$2$};
  \node[squarednode](2) at (2.6,0){$4$};
    \node[squarednode](3) at (3.9,0){$7$};
     \node at (2.6,-0.7){$C_{\{2,4,7\}}$};
    \draw[thick,double,->](1)--(2);
       \draw[ thick,double,->](2)--(3);
        \draw[ thick,double,->](3)--(3.9,0.7)--(1.3,0.7)--(1);
 \end{tikzpicture} \quad 
   \begin{tikzpicture}[
roundnode/.style={circle, draw=black!60, fill=gray!10, very thick, minimum size=7mm},
squarednode/.style={rectangle, draw=black!60, fill=gray!5, very thick, minimum size=5mm},]
 \node[squarednode](1) at (1.3,0){$2$};
  \node[squarednode](2) at (2.6,0){$4$};
    \node[squarednode](3) at (3.9,0){$7$};
     \node at (2.6,-0.7){$P_{\{2,4,7\}}$};
    \draw[thick,->](1)--(2);
       \draw[ thick,->](2)--(3);
 \end{tikzpicture} \quad
  \begin{tikzpicture}[
roundnode/.style={circle, draw=black!60, fill=gray!10, very thick, minimum size=7mm},
squarednode/.style={rectangle, draw=black!60, fill=gray!5, very thick, minimum size=5mm},]
 \node[squarednode](1) at (1.3,0){$2$};
  \node[squarednode](2) at (2.6,0){$4$};
    \node[squarednode](3) at (3.9,0){$7$};
    \node at (2.6,-0.7){$Q_{\{2,4,7\}}$};
    \draw[thick,double,->](1)--(2);
       \draw[ thick,double,->](2)--(3);
 \end{tikzpicture} \quad
  \begin{tikzpicture}[
roundnode/.style={circle, draw=black!60, fill=gray!10, very thick, minimum size=7mm},
squarednode/.style={rectangle, draw=black!60, fill=gray!5, very thick, minimum size=5mm},]
 \node[squarednode](1) at (1.3,0){$2$};
  \node[squarednode](2) at (2.6,0){$4$};
    \node[squarednode](3) at (3.9,0){$7$};
     \node at (2.6,-0.7){$K_{\{2,4,7\}}$};
    \draw[thick,dashed,->](1)--(2);
       \draw[ thick,dashed,->](2)--(3);
        \draw[ thick,dashed,->](1)--(1.3,0.7)--(3.9,0.7)--(3);
 \end{tikzpicture}
  \caption{\fathree{The labelled edge-coloured digraphs $C_{\{2,4,7\}},P_{\{2,4,7\}},Q_{\{2,4,7\}}, K_{\{2,4,7\}}$}}\label{graphsn}
  \end{center}
  \end{figure}
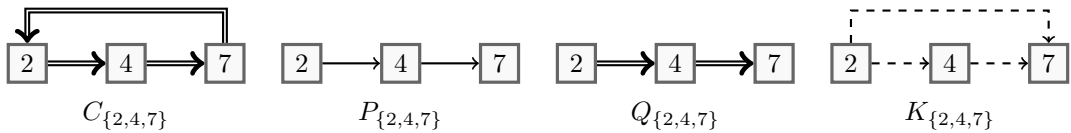
  
\fa{
In the rest of the paper, if $\kappa$ is any type of colouring of a labelled graph, digraph, or edge-coloured digraph with vertex set $[n]$, we define 
\begin{equation}\label{eq:xkappan}
\bx_\kappa= \bx_{\kappa(1)}\bx_{\kappa(2)} \cdots \bx_{\kappa(n)}.
\end{equation} 
}

\begin{definition} {(The generalized chromatic function of a labelled edge-coloured digraph in noncommuting variables)}\label{GCFn}
Let $\bG$ be a labelled edge-coloured digraph with vertex set $[n]$. Define
$$
\GY_{\bG}({\bf x},t)=\sum_{\kappa\in \mathcal{C}(\bG)}t^{{\rm asc}(\kappa)} {\bf x}_{\kappa} 
$$
 \fa{where $\bx_\kappa$ is defined by Equation \ref{eq:xkappan}} and   $${\rm asc}(\kappa)=|\{(a,b)\in E(\bG):  \kappa(a)<\kappa(b)\}|.$$
The \fss{power series} $$\GY_{\bG}({\bf x})=\GY_{\bG}({\bf x},1)$$ is called the \emph{generalized chromatic function} of $\bG$.  
\end{definition}

Consider the following commutation map, 
$$
\begin{array}{cccc}
\rho :& \mathbb{Q}\langle \langle \bx_1,\bx_2,\ldots \rangle \rangle&  \rightarrow &   \mathbb{Q}[[ x_1,x_2,\ldots ]]\\
& \bx_i &\mapsto & x_i.
\end{array}
$$

By definition we have the following.
\begin{proposition}\label{prop:commute-map}
For any labelled edge-coloured digraph $\bG$, we have 
$$\rho (\GY_{\bG}({\bf x},t))=\GX_{G}(x,t),$$ where $G$ is $\bG$ with labels removed. 
\end{proposition}

\paragraph{Gebhard-Sagan's chromatic symmetric function in noncommuting variables \cite{GS}:}
Let $H$ be a labelled graph with vertex set $[n]$.  Define 
$$Y_{H}(\bx)=\sum_{\kappa \in \mathcal{C}(H)} {\bf x}_\kappa,$$ 
 \fa{where $\bx_\kappa$ is defined by Equation \ref{eq:xkappan}}.

By definition we have the following.
\begin{proposition} \label{prop:GSasGCF}
Let $H$ be a labelled graph with vertex set $[n]$. Then 
$$Y_H(\bx)=\GY_{\bG}(\bx)$$ where $\bG$ is any labelled edge-coloured digraph that has only dashed edges and mapping the vertex $a$ of the underlying graph of $\bG$ to the vertex $a$ of $H$  gives an isomorphism.
 \end{proposition}

Let $A$ be a nonempty subset of $\mathbb{P}$. Define $\mathfrak{S}_A$ to be the set of all bijections from $A$ to itself. Let $\bG$ be a labelled edge-coloured digraph with vertex set $A$. Then for \fathree{$\sigma\in \mathfrak{S}_A$}, define 
$\sigma\circ {\bG}$ to be the labelled edge-coloured digraph obtained by replacing each vertex $a$ of $\bG$ by $\sigma(a)$. The \emph{symmetrized generalized chromatic function} of ${\bG}$ is 
$$\mathfrak{S}\GY_{\bG}(\bx)=\sum_{\sigma\in \mathfrak{S}_A}\GY_{\sigma\circ \bG} (\bx).$$
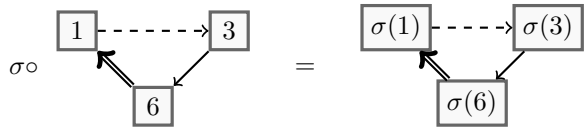
\begin{figure}[H]
\begin{center} 
 \begin{tikzpicture}[
roundnode/.style={circle, draw=black!60, fill=gray!10, very thick, minimum size=7mm},
squarednode/.style={rectangle, draw=black!60, fill=gray!5, very thick, minimum size=5mm},]
\node at (0.3,0.5){$\sigma \circ$};
 \node[squarednode](1) at (1,1){$1$};
  \node[squarednode](3) at (2,0){$6$};
    \node[squarednode](2) at (3,1){$3$};
    \draw[thick, dashed,->](1)--(2);
       \draw[ thick,->](2)--(3);
        \draw[ thick,double,->](3)--(1);
        \node at (4,0.5){$\quad = \quad $};
 \end{tikzpicture}
  \begin{tikzpicture}[
roundnode/.style={circle, draw=black!60, fill=gray!10, very thick, minimum size=7mm},
squarednode/.style={rectangle, draw=black!60, fill=gray!5, very thick, minimum size=5mm},]
 \node[squarednode](1) at (1,1){$\sigma(1)$};
  \node[squarednode](3) at (2,0){$\sigma(6)$};
    \node[squarednode](2) at (3,1){$\sigma(3)$};
    \draw[thick, dashed,->](1)--(2);
       \draw[ thick,->](2)--(3);
        \draw[ thick,double,->](3)--(1);
 \end{tikzpicture}
  \caption{Action of a permutation $\sigma\in \mathfrak{S}_{\{1,3,6\}}$ on the  labelled edge-coloured digraph in Figure \ref{ex2}}\label{ex:Action}
 \end{center} 
 \end{figure}

\section{Product and Coproduct Formulas for Generalized Chromatic Functions in Noncommuting Variables}\label{sec:prodsfornc}

We now establish the product and coproduct formulas for generalized chromatic functions in noncommuting variables. 

Let $\bG$ be a labelled edge-coloured digraph, and let $n$ be a positive integer. Then $\bG+n$ is the labelled edge-coloured digraph obtained by replacing each vertex $a$ of $\bG$ by $a+n$. 
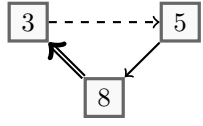
\begin{figure}[H]
\begin{center} 
 \begin{tikzpicture}[
roundnode/.style={circle, draw=black!60, fill=gray!10, very thick, minimum size=7mm},
squarednode/.style={rectangle, draw=black!60, fill=gray!5, very thick, minimum size=5mm},]
 \node[squarednode](1) at (1,1){$3$};
  \node[squarednode](3) at (2,0){$8$};
    \node[squarednode](2) at (3,1){$5$};
    \draw[thick, dashed,->](1)--(2);
       \draw[ thick,->](2)--(3);
        \draw[ thick,double,->](3)--(1);
 \end{tikzpicture}
  \caption{$\bG+2$ where $\bG$ is the  labelled edge-coloured digraph in Figure \ref{ex2}}\label{ex:Shift}
 \end{center} 
\end{figure}

\begin{proposition}\label{prodn}
Let $\bG_1$ and $\bG_2$ be labelled edge-coloured digraphs.  Then 
$$
\GY_{\bG_1 }({\bf x},t)\GY_{\bG_2}({\bf x},t)=\GY_{\bG_1\uplus  (\bG_2+|V(\bG_1)|) }({\bf x},t).
$$
\end{proposition} 

Let $\bG$ be a labelled edge-coloured digraph. The \emph{standardization} of $\bG$, denoted ${\rm std}(\bG)$, is the labelled edge-coloured digraph obtained by replacing the $i$th smallest element of $V(\bG)$ by $i$.
\begin{figure}[H]
\begin{center} 
 \begin{tikzpicture}[
roundnode/.style={circle, draw=black!60, fill=gray!10, very thick, minimum size=7mm},
squarednode/.style={rectangle, draw=black!60, fill=gray!5, very thick, minimum size=5mm},]
 \node[squarednode](1) at (1,1){$1$};
  \node[squarednode](3) at (2,0){$3$};
    \node[squarednode](2) at (3,1){$2$};
    \draw[thick, dashed,->](1)--(2);
       \draw[ thick,->](2)--(3);
        \draw[ thick,double,->](3)--(1);
 \end{tikzpicture}
  \caption{${\rm std}(\bG)$ where $\bG$ is the  labelled edge-coloured digraph in Figure \ref{ex2}}\label{ex:Shift}
 \end{center} 
 \end{figure}

  By a proof analogous to that of Theorem \ref{coprod}, we have the following coproduct formula for generalized chromatic functions of labelled edge-coloured digraphs.

\begin{proposition}\label{coprodn}
Let $\bG$ be a labelled edge-coloured digraph. Then
$$\Delta(\GY_{\bG}({\bf x}))=\sum  \GY_{{\rm std}(\bG|_{V(\bG)-V(\bF)})}({\bf x})\otimes \GY_{{\rm std}(\bF)}({\bf x}),$$ where the sum runs over all $\{\rightarrow,\Rightarrow\}$-induced subdigraphs $\bF$ of $\bG$.
\end{proposition} 

\section{Bases for $\NCSym(\bx)$ and $\NCQSym(\bx)$  Using Generalized Chromatic Functions in Noncommuting Variables}\label{basesn}

In this section, we show that different bases for the Hopf algebras $\NCSym(\bx)$ and \fathree{$\NCQSym(\bx)$}  can be realized as generalized chromatic functions in noncommuting variables of certain labelled edge-coloured digraphs. Since these bases are less well-known than their \fathree{commutative} counterparts, we provide their classical definitions taken from \cite{RS} for $\NCSym(\bx)$ and \cite{BZ, UpFun} for $\NCQSym(\bx)$.

\paragraph{Monomial basis $\{{\bf m}_\Pi\}$ of $\NCSym(\bx)$.} 
Given $\Pi=\Pi_1/\Pi_2/\cdots/\Pi_l\vdash [n]$, define the 
\emph{monomial symmetric function} in noncommuting variables to be
$$
{\bf m}_{\Pi}
= \sum_{(i_{1},i_{2},\ldots ,i_{n})} \bx_{i_{1}}\bx_{i_{2}}\cdots \bx_{i_{n}},
$$
where the sum is over all $n$-tuples   $(i_1,i_2,\ldots,i_n)$ with
$i_j=i_k$ if and only if $j$ and $k$ are in the same
  block of $\Pi$. 
For example, 
$$
{\bf m}_{13/24}=\bx_1\bx_2\bx_1\bx_2 +\bx_2\bx_1\bx_2\bx_1  +\bx_1\bx_3\bx_1\bx_3+ 
\bx_3\bx_1\bx_3\bx_1+\bx_2\bx_3\bx_2\bx_3+\bx_3\bx_2\bx_3\bx_2  +\cdots.
$$ 

By definition we have  that ${\bf m}_{\Pi}=\GY_{\bG}(\bx)$ where  
$$
\begin{tikzpicture} 
\node (a) at (-0.2,0){${\bG}=$};
\node (b) at (1,0){$C_{\Pi_1}$};
\node (c) at (1.7,0){$\arx$};
\node (d) at (2.5,0){$C_{\Pi_2}$};
\node (e) at (3.2,0){$\arx$};
\node (f) at (4,0){$\dots$};
\node (g) at (4.7,0){$\arx$};
\node (h) at (5.6,0){$C_{\Pi_l}.$};
\end{tikzpicture} 
$$For example, if $\bG$ is the following labelled  edge-coloured digraph 
$$
  
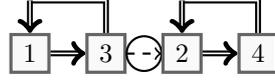
\begin{tikzpicture}[
roundnode/.style={circle, draw=black!60, fill=gray!10, very thick, minimum size=7mm},
squarednode/.style={rectangle, draw=black!60, fill=gray!5, very thick, minimum size=5mm},]
      \node[squarednode](2) at (2,0){$1$};
        \node[squarednode](3) at (3,0){$3$};
               \node(a) at (3.5,0){$\arx$};
          \node[squarednode](4) at (4,0){$2$};
                    \node[squarednode](5) at (5,0){$4$};
\draw[thick,double,->](2)--(3);
\draw[thick,double](3)--(3,0.7)--(2,0.7);
\draw[thick,double,->](2,0.7)--(2);
\draw[thick,double,->](4)--(5);
\draw[thick,double](5)--(5,0.7)--(4,0.7);
\draw[thick,double,->](4,0.7)--(4);

%
  \end{tikzpicture}
$$
  then $\GY_{\bG}(\bx)=\mathbf{m}_{13/24}$.

\paragraph{Power sum basis $\{{\bf p}_\Pi\}$ of $\NCSym(\bx)$.} Given $\Pi=\Pi_1/\Pi_2/\cdots/\Pi_l\vdash [n]$, define the \emph{power sum symmetric function} in noncommuting variables to be
$${\bf p}_{\Pi}=\sum_{(i_1,i_2,\ldots,i_n)}\bx_{i_1}\bx_{i_2}\cdots \bx_{i_n},$$ 
where $i_j=i_k$ if $j,k$ are in the same block of $\Pi$.
For example, 
$$\mathbf{p}_{13/24}=\bx_1\bx_2\bx_1\bx_2+\bx_2\bx_1\bx_2\bx_1+\bx_1^4+\bx_2^4+\cdots.$$
By definition we have that ${\bf p}_{\Pi}=\GY_{\bG}(\bx)$ where 
$$
\begin{tikzpicture} 
\node (a) at (-0.2,0){${\bG}=$};
\node (b) at (1,0){$C_{\Pi_1}$};
\node (c) at (1.7,0){$\biguplus$};
\node (d) at (2.5,0){$C_{\Pi_2}$};
\node (e) at (3.2,0){$\biguplus$};
\node (f) at (4,0){$\dots$};
\node (g) at (4.7,0){$\biguplus$};
\node (h) at (5.6,0){$C_{\Pi_l}.$};
\end{tikzpicture} 
$$
For example, if $\bG$ is the following labelled edge-coloured digraph
$$
  
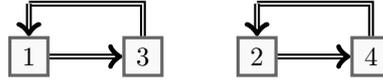
\begin{tikzpicture}[
roundnode/.style={circle, draw=black!60, fill=gray!10, very thick, minimum size=7mm},
squarednode/.style={rectangle, draw=black!60, fill=gray!5, very thick, minimum size=5mm},]
 \node[squarednode](1) at (1,0){1};
  \node[squarednode](2) at (2.5,0){3};
    \node[squarednode](3) at (4,0){2};
       \node[squarednode](4) at (5.5,0){4};
    \draw[thick,double,->](1)--(2);
       \draw[thick,double](2)--(2.5,0.7)--(1,0.7);
              \draw[thick,double,->](1,.7)--(1);
                  \draw[thick,double,->](3)--(4);
       \draw[thick,double](4)--(5.5,0.7)--(4,0.7);
              \draw[thick,double,->](4,.7)--(3);
 \end{tikzpicture}
$$
then $\GY_{\bG}(\bx)=\mathbf{p}_{13/24}$.

 \begin{remark}
  Note that the underlying edge-coloured digraphs for the augmented monomial basis elements of $\Sym(x)$ and the monomial basis elements of $\NCSym(\bx)$ are the same.  
Moreover,  the underlying edge-coloured digraphs for power sum basis elements of $\Sym(x)$ and power sum basis elements in $\NCSym(\bx)$ are the same. Thus by Proposition \ref{prop:commute-map}, we recover \cite[Theorem 2.1 (i), (ii), (iii)]{RS}, where the third part is recovered from the first interpretation below.
\end{remark}

\paragraph{Elementary basis $\{{\bf e}_{\Pi}\}$ of $\NCSym(\bx)$.} Given $\Pi=\Pi_1/\Pi_2/\cdots/\Pi_l\vdash [n]$, define 
the  \emph{elementary symmetric function} in noncommuting
variables to be
$${\bf e}_{\Pi}= \sum_{(i_{1},i_{2},\ldots ,i_{n})} \bx_{i_{1}}\bx_{i_{2}}\cdots \bx_{i_{n}},$$
where $i_j\neq i_k$ if $j,k$ are in the same block of $\Pi$.
For example,
$$
{\bf e}_{13/24}=\bx_1\bx_1\bx_2\bx_2+\bx_2\bx_2\bx_1\bx_1+\bx_1\bx_2\bx_2\bx_1+\fa{\bx_2\bx_1\bx_1\bx_2}+\bx_1\bx_2\bx_3\bx_4+\cdots.
$$

It follows by definition that the elementary symmetric functions in noncommuting variables can be written in two ways using labelled edge-coloured digraphs. 
The first is that
${\bf e}_{\Pi}=\GY_{\bG}(\bx)$ where 
$$
\begin{tikzpicture} 
\node (a) at (-0.2,0){${\bG}=$};
\node (b) at (1,0){$K_{\Pi_1}$};
\node (c) at (1.7,0){$\biguplus$};
\node (d) at (2.5,0){$K_{\Pi_2}$};
\node (e) at (3.2,0){$\biguplus$};
\node (f) at (4,0){$\dots$};
\node (g) at (4.7,0){$\biguplus$};
\node (h) at (5.6,0){$K_{\Pi_l}.$};
\end{tikzpicture} 
$$
For example, if $\bG$ is the following labelled edge-coloured digraph
$$
  \begin{tikzpicture}[
roundnode/.style={circle, draw=black!60, fill=gray!10, very thick, minimum size=7mm},
squarednode/.style={rectangle, draw=black!60, fill=gray!5, very thick, minimum size=5mm},]
 \node[squarednode](1) at (1,1){1};
  \node[squarednode](2) at (2.5,1){3};
    \node[squarednode](3) at (3.5,1){2};
       \node[squarednode](4) at (5,1){4};
    \draw[thick,dashed,->](1)--(2);
                  \draw[thick,dashed,->](3)--(4);
 \end{tikzpicture}
 $$
then $\GY_{\bG}(\bx)=\mathbf{e}_{13/24}$.

The second is that
$$\mathbf{e}_{\Pi}=\sum_{(\sigma_1,\ldots,\sigma_l)\in \mathfrak{S}_{\Pi_1}\times \cdots \times \mathfrak{S}_{\Pi_l}} \GY_{(\uplus_{i=1}^l \sigma_i\circ P_{\Pi_i})}(\bx).$$

 For example, $\mathbf{e}_{13/24}$ is equal to
$$ 
\GY_{ 
  \begin{tikzpicture}[
roundnode/.style={circle, draw=black!60, fill=gray!10, very thick, minimum size=7mm},
squarednode/.style={rectangle, draw=black!60, fill=gray!5, very thick, minimum size=5mm},]
 \node[squarednode](1) at (1,1){1};
  \node[squarednode](2) at (2.5,1){3};
    \node[squarednode](3) at (3.5,1){2};
       \node[squarednode](4) at (5,1){4};
    \draw[thick,->](1)--(2);
                  \draw[ thick,->](3)--(4);
 \end{tikzpicture}
 }({\bf x})+\GY_{ 
  \begin{tikzpicture}[
roundnode/.style={circle, draw=black!60, fill=gray!10, very thick, minimum size=7mm},
squarednode/.style={rectangle, draw=black!60, fill=gray!5, very thick, minimum size=5mm},]
 \node[squarednode](1) at (1,1){3};
  \node[squarednode](2) at (2.5,1){1};
    \node[squarednode](3) at (3.5,1){2};
       \node[squarednode](4) at (5,1){4};
    \draw[thick,->](1)--(2);
                  \draw[thick,->](3)--(4);
 \end{tikzpicture}
 }({\bf x})$$
$$
+\GY_{ 
  \begin{tikzpicture}[
roundnode/.style={circle, draw=black!60, fill=gray!10, very thick, minimum size=7mm},
squarednode/.style={rectangle, draw=black!60, fill=gray!5, very thick, minimum size=5mm},]
 \node[squarednode](1) at (1,1){1};
  \node[squarednode](2) at (2.5,1){3};
    \node[squarednode](3) at (3.5,1){4};
       \node[squarednode](4) at (5,1){2};
    \draw[thick,->](1)--(2);
                  \draw[thick,->](3)--(4);
 \end{tikzpicture}
 }({\bf x})
 +\GY_{ 
  \begin{tikzpicture}[
roundnode/.style={circle, draw=black!60, fill=gray!10, very thick, minimum size=7mm},
squarednode/.style={rectangle, draw=black!60, fill=gray!5, very thick, minimum size=5mm},]
 \node[squarednode](1) at (1,1){3};
  \node[squarednode](2) at (2.5,1){1};
    \node[squarednode](3) at (3.5,1){4};
       \node[squarednode](4) at (5,1){2};
    \draw[thick,->](1)--(2);
                  \draw[thick,->](3)--(4);
 \end{tikzpicture}
 }({\bf x}).
$$

\paragraph{Complete homogeneous basis $\{{\bf h}_\Pi\}$ of $\NCSym(\bx)$.}
 For set partitions $\Pi$ and $\Omega$ of $[n]$, let $\Omega\leq \Pi$ if each block of $\Omega$ is contained in some block of $\Pi$. The set of all set partitions of $[n]$ with this partial ordering gives a lattice; the meet (greatest lower bound) and join (least upper bound) operations of this lattice are denoted by $\wedge$ 
and $\vee$, respectively. 
The \emph{complete homogeneous symmetric function} in noncommuting variables is
$$
{\bf h}_\Pi=\sum_{\Omega \vdash [n]} \fathree{(\lambda(\Omega \wedge \Pi))!} {\bf m}_\Omega.
$$
For example, 
$$
{\bf h}_{13/24}={\bf m}_{1/2/3/4}+{\bf m}_{12/3/4}+2{\bf m}_{13/2/4}+{\bf m}_{14/2/3}
			+{\bf m}_{1/23/4}+2{\bf m}_{1/24/3}+{\bf m}_{1/2/34}$$
			$$
 +{\bf m}_{12/34}+4{\bf m}_{13/24}+{\bf m}_{14/23}
+2{\bf m}_{123/4}+2{\bf m}_{124/3}+2{\bf m}_{134/2}+2{\bf m}_{1/234}+4{\bf m}_{1234}.
$$

Using labelled edge-coloured digraphs we now present the complete homogeneous basis. Let $\Pi=\Pi_1/\Pi_2/\cdots/\Pi_l$ be a set partition of $[n]$.
By \cite[Lemma 2.14]{ALvW}, we have that
$$\mathbf{h}_{\Pi}=\sum_{(\sigma_1,\ldots,\sigma_l)\in \mathfrak{S}_{\Pi_1}\times \cdots \times \mathfrak{S}_{\Pi_l}} \GY_{(\uplus_{i=1}^l \sigma_i\circ Q_{\Pi_i})}(\bx).$$
 For example, $\mathbf{h}_{13/24}$ is equal to
$$ 
\GY_{ 
  \begin{tikzpicture}[
roundnode/.style={circle, draw=black!60, fill=gray!10, very thick, minimum size=7mm},
squarednode/.style={rectangle, draw=black!60, fill=gray!5, very thick, minimum size=5mm},]
 \node[squarednode](1) at (1,1){1};
  \node[squarednode](2) at (2.5,1){3};
    \node[squarednode](3) at (3.5,1){2};
       \node[squarednode](4) at (5,1){4};
    \draw[ thick,double,->](1)--(2);
                  \draw[thick,double,->](3)--(4);
 \end{tikzpicture}
 }({\bf x})+\GY_{ 
  \begin{tikzpicture}[
roundnode/.style={circle, draw=black!60, fill=gray!10, very thick, minimum size=7mm},
squarednode/.style={rectangle, draw=black!60, fill=gray!5, very thick, minimum size=5mm},]
 \node[squarednode](1) at (1,1){3};
  \node[squarednode](2) at (2.5,1){1};
    \node[squarednode](3) at (3.5,1){2};
       \node[squarednode](4) at (5,1){4};
    \draw[thick,double,->](1)--(2);
                  \draw[thick,double,->](3)--(4);
                   \end{tikzpicture}
 }({\bf x})$$
$$
+\GY_{ 
  \begin{tikzpicture}[
roundnode/.style={circle, draw=black!60, fill=gray!10, very thick, minimum size=7mm},
squarednode/.style={rectangle, draw=black!60, fill=gray!5, very thick, minimum size=5mm},]
 \node[squarednode](1) at (1,1){1};
  \node[squarednode](2) at (2.5,1){3};
    \node[squarednode](3) at (3.5,1){4};
       \node[squarednode](4) at (5,1){2};
    \draw[thick,double,->](1)--(2);
                  \draw[thick,double,->](3)--(4);
                   \end{tikzpicture}
 }({\bf x})
 +\GY_{ 
  \begin{tikzpicture}[
roundnode/.style={circle, draw=black!60, fill=gray!10, very thick, minimum size=7mm},
squarednode/.style={rectangle, draw=black!60, fill=gray!5, very thick, minimum size=5mm},]
 \node[squarednode](1) at (1,1){3};
  \node[squarednode](2) at (2.5,1){1};
    \node[squarednode](3) at (3.5,1){4};
       \node[squarednode](4) at (5,1){2};
    \draw[thick,double,->](1)--(2);
                  \draw[thick,double,->](3)--(4);
 \end{tikzpicture}
 }({\bf x}).
$$

\fathree{
\begin{remark} Given a partition $\lambda$, the generalized chromatic functions of $P_{\lambda_i}$'s give the elementary symmetric function in commuting variables, $e_\lambda$. Likewise, given a set partition $\Pi$, the symmetrized generalized chromatic functions of $P_{\Pi_i}$'s give the elementary symmetric function in noncommuting variables ${\bf e}_\Pi$. Similarly, by looking at $Q_{\lambda_i}$'s and $Q_{\Pi_i}$'s we can obtain the complete homogeneous symmetric functions $h_\lambda$ and ${\bf h}_\Pi$. 
\end{remark}
}

\paragraph{Rosas-Sagan Schur functions of $\NCSym(\bx)$.} Rosas and Sagan in \cite{RS}, as an analogy for the monomial, power sum, elementary and homogeneous bases of $\Sym(x)$, introduced the above bases for $\NCSym(\bx)$, recalling the elementary basis from the work of Wolf \cite{Wolf}. Their proposed analogy for Schur functions did not produce enough distinct elements to make a basis for 
$\NCSym(\bx)$. However, their functions have a natural realization in terms of generalized chromatic functions, so we include them here.

Let  $\Pi=\Pi_1/\Pi_2/\cdots/\Pi_l$ be a set partition of $[n]$ with $\lambda=\lambda(\Pi)$. Let $\bG_\lambda$ denote the labelled edge-coloured digraph obtained by replacing the vertex $(i,j)$ in $G_\lambda$ by $\lambda_1+\cdots+\lambda_{i-1}+j$.
Now define
$${\bf S}_{\Pi}=\sum_{\sigma\in \mathfrak{S}_n} \GY_{\sigma\circ \bG_\lambda}({\bf x}).$$
Then ${\bf S}_{\Pi} \in \NCSym(\bx)$, and ${\bf S}_{\Pi}={\bf S}_{\Omega}$ if and only if $\lambda(\Pi)=\lambda(\Omega)=\lambda$. Note that ${\bf S}_{\Pi}$ is the same as \svwt{${S}_\lambda$} introduced in \cite{RS}. Moreover, $\rho({\bf S}_{\Pi})=n! s_\lambda$. However, the set $$\{{\bf S}_{\Pi}: \Pi\vdash [n]\}$$  is not a basis for $\NCSym_n(\bx)$ since the dimension of the space spanned by this set is equal to the number of partitions, which is less than the number of set partitions for $n>2$.

\paragraph{Monomial basis $\{{\bf M}_\Phi\}$ of $\NCQSym(\bx)$.}
Given {$\Phi=(\Phi_1|\Phi_2|\cdots|\Phi_k)\vDash [n]$,  define the \emph{monomial quasisymmetric function} in noncommuting variables  to be
$${\bf M}_{\Phi}=\sum_{(i_1,i_2,\ldots,i_n)} \bx_{i_1}\bx_{i_2}\cdots \bx_{i_n}$$ where the sum runs over all tuples $(i_1,i_2,\ldots,i_n)$ such that: 
\sxl{\begin{itemize}
\item[(i)] $i_j=i_\ell$ if $j$ and $\ell$ are  in the same block of $\Phi$.
\item[(ii)] $i_j<i_\ell$ if \fss{$j\in \Phi_p$ and $\ell\in \Phi_q$ with $p<q$.}
\end{itemize}}

For example,
$${\bf M}_{(13|24)}= \bx _1 \bx _2 \bx _1 \bx _2 + \bx _1 \bx _3 \bx _1 \bx _3 + \bx _2 \bx _3 \bx _2 \bx _3+\cdots . $$

By definition, we have that ${\bf M}_{\Phi}=\GY_{\bG}(\bx),$ where  
 $$
\begin{tikzpicture} 
\node (a) at (-0.2,0){${\bG}=$};
\node (b) at (1,0){$C_{\Phi_1}$};
\node (c) at (1.7,0){$\ary$};
\node (d) at (2.5,0){$C_{\Phi_2}$};
\node (e) at (3.2,0){$\ary$};
\node (f) at (4,0){$\dots$};
\node (g) at (4.7,0){$\ary$};
\node (h) at (5.6,0){$C_{\Phi_k}.$};
\end{tikzpicture} 
$$

For example, if $\bG$ is the following labelled edge-coloured digraph 
  $$
  
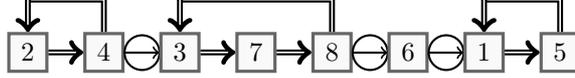
\begin{tikzpicture}[
roundnode/.style={circle, draw=black!60, fill=gray!10, very thick, minimum size=7mm},
squarednode/.style={rectangle, draw=black!60, fill=gray!5, very thick, minimum size=5mm},]
 \node[squarednode](1) at (0,0){$2$};
  \node[squarednode](2) at (1,0){$4$};
  \node(a) at (1.5,0){$\ary$};
    \node[squarednode](3) at (2,0){$3$};
      \node[squarednode](4) at (3,0){$7$};
        \node[squarednode](5) at (4,0){$8$};
          \node(a) at (4.5,0){$\ary$};
          \node[squarednode](6) at (5,0){$6$};
            \node(a) at (5.5,0){$\ary$};
                    \node[squarednode](7) at (6,0){$1$};
                                        \node[squarednode](8) at (7,0){$5$};
\draw[thick,double,->](1)--(2);
\draw[thick,double](2)--(1,0.7)--(0,0.7);
\draw[thick,double,->](0,0.7)--(1);
\draw[thick,double](5)--(4,0.7)--(2,0.7);
\draw[thick,double,->](2,0.7)--(3);
\draw[thick,double](8)--(7,0.7)--(6,0.7);
\draw[thick,double,->](6,0.7)--(7);

\draw[thick,double,->](3)--(4);
\draw[thick,double,->](4)--(5);
\draw[thick,double,->](7)--(8);
  \end{tikzpicture}
  $$then $\GY_{\bG}(x)={\bf M}_{(24|378|6|15)}$.
  
  To the best of our knowledge the following basis is not in the literature, however, similar bases exist \cite{BZ,UpFun}. Therefore we will define the basis elements in terms of generalized chromatic functions in noncommuting variables first, before deriving an explicit formula for them, and establishing that they are a basis for $\NCQSym (\bx)$. Since they can be defined using generalized chromatic functions in noncommuting variables their product and coproduct formulas follow from Propositions~\ref{prodn} and \ref{coprodn}.
  
\paragraph{Fundamental basis  $\{{\bf F}_\Phi\}$ of $\NCQSym(\bx)$.}
Given $\Phi=(\Phi_1|\Phi_2|\cdots|\Phi_k)\vDash [n]$,   define the \emph{fundamental quasisymmetric function} in noncommuting variables to be
 $$\bF_{\Phi}=\GY_{\bG}(\bx),$$ where 
  $$
\begin{tikzpicture} 
\node (a) at (-0.2,0){${\bG}=$};
\node (b) at (1,0){$Q_{\Phi_1}$};
\node (c) at (1.7,0){$\ary$};
\node (d) at (2.5,0){$Q_{\Phi_2}$};
\node (e) at (3.2,0){$\ary$};
\node (f) at (4,0){$\dots$};
\node (g) at (4.7,0){$\ary$};
\node (h) at (5.6,0){$Q_{\Phi_k}.$};
\end{tikzpicture} 
$$
For example, if $\Phi = (24|378|6|15)$ then $\bG$ is the following labelled edge-coloured digraph
$$
  \begin{tikzpicture}[
roundnode/.style={circle, draw=black!60, fill=gray!10, very thick, minimum size=7mm},
squarednode/.style={rectangle, draw=black!60, fill=gray!5, very thick, minimum size=5mm},]
 \node[squarednode](1) at (0,0){$2$};
  \node[squarednode](2) at (1,0){$4$};
  \node(a) at (1.5,0){$\ary$};
    \node[squarednode](3) at (2,0){$3$};
      \node[squarednode](4) at (3,0){$7$};
        \node[squarednode](5) at (4,0){$8$};
          \node(a) at (4.5,0){$\ary$};
          \node[squarednode](6) at (5,0){$6$};
            \node(a) at (5.5,0){$\ary$};
                    \node[squarednode](7) at (6,0){$1$};
                                        \node[squarednode](8) at (7,0){$5$};
\draw[thick,double,->](1)--(2);

\draw[thick,double,->](3)--(4);
\draw[thick,double,->](4)--(5);
\draw[thick,double,->](7)--(8);
  \end{tikzpicture}
  $$
and $\GY_{\bG}(\bx)=\bF_{(24|378|6|15)}.$

By comparing the labelled edge-coloured digraphs in the realizations of $\bF _\Phi$ and ${\bf M}_\Phi$ as generalized chromatic functions in noncommuting variables, we immediately get the following.

\begin{proposition}\label{prop:NCFasM} Let $\Phi$ be a set composition. Then
$$\bF_\Phi=\sum_{\Psi \text{\svwt{~reforms~}} \Phi}  {\bf M}_\Psi$$and hence $\{\bF_\Phi\}$ is a basis for $\NCQSym(\bx).$ \end{proposition}

For example,   $$\bF_{(13|24)} = {\bf M} _{(13|24)}+{\bf M} _{(1|3|24)}+{\bf M} _{(13|2|4)}+{\bf M} _{(1|3|2|4)}.$$

\paragraph{Upper-fundamental basis $\{\overline{\bF}_{\Phi}\}$ of $\NCQSym(\bx)$.} This basis also appears in \cite{UpFun} as the $L$ basis, but we \fss{reinterpret} it here as generalized chromatic functions. Given $\Phi=(\Phi_1|\Phi_2|\cdots|\Phi_k)\vDash[n]$,  define the \emph{upper-fundamental quasisymmetric function} in noncommuting variables to be
$$\overline{\bF}_\Phi=\sum_{\Psi \text{\svwt{~corrupts~}} \Phi} \sxl{\bf{M}}_\Psi.$$
For example,
$$\fss{\overline{\bF}_{(13|24)} = {\bf M} _{(13|24)}+{\bf M} _{(1324)}= {\bf M} _{(13|24)}+{\bf M} _{(1234)}.}$$

By definition, we have that $\overline{\bF}_{\Phi}=\GY_{\bG}(\bx),$ where 
$$
\begin{tikzpicture} 
\node (a) at (-0.2,0){${\bG}=$};
\node (b) at (1,0){$C_{\Phi_1}$};
\node (c) at (1.7,0){$\arz$};
\node (d) at (2.5,0){$C_{\Phi_2}$};
\node (e) at (3.2,0){$\arz$};
\node (f) at (4,0){$\dots$};
\node (g) at (4.7,0){$\arz$};
\node (h) at (5.6,0){$C_{\Phi_k}.$};
\end{tikzpicture} 
$$
For example, if $\bG$ is the following labelled edge-coloured digraph
$$
  
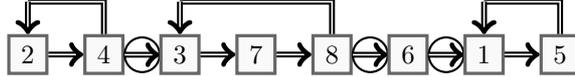
\begin{tikzpicture}[
roundnode/.style={circle, draw=black!60, fill=gray!10, very thick, minimum size=7mm},
squarednode/.style={rectangle, draw=black!60, fill=gray!5, very thick, minimum size=5mm},]
 \node[squarednode](1) at (0,0){$2$};
  \node[squarednode](2) at (1,0){$4$};
  \node(a) at (1.5,0){$\arz$};
    \node[squarednode](3) at (2,0){$3$};
      \node[squarednode](4) at (3,0){$7$};
        \node[squarednode](5) at (4,0){$8$};
          \node(a) at (4.5,0){$\arz$};
          \node[squarednode](6) at (5,0){$6$};
            \node(a) at (5.5,0){$\arz$};
                    \node[squarednode](7) at (6,0){$1$};
                                        \node[squarednode](8) at (7,0){$5$};
\draw[thick,double,->](1)--(2);
\draw[thick,double](2)--(1,0.7)--(0,0.7);
\draw[thick,double,->](0,0.7)--(1);
\draw[thick,double](5)--(4,0.7)--(2,0.7);
\draw[thick,double,->](2,0.7)--(3);
\draw[thick,double](8)--(7,0.7)--(6,0.7);
\draw[thick,double,->](6,0.7)--(7);

\draw[thick,double,->](3)--(4);
\draw[thick,double,->](4)--(5);
\draw[thick,double,->](7)--(8);
  \end{tikzpicture}
  $$
then $\GY_{\bG}(\bx)=\overline{\bF}_{(24|378|6|15)}.$

\begin{remark} 
Note that the underlying edge-coloured digraphs of monomial (fundamental and upper-fundamental, respectively) bases of $\QSym(x)$ and $\NCQSym(\bx)$ are the same. Moreover, for any set composition $\Phi$, 
$$
\rho(\bM_\Phi)=M_{\alpha(\Phi)},\quad \quad \rho(\bF_\Phi)=F_{\alpha(\Phi)},~~ \text{and} \quad \quad \rho(\overline{\bF}_\Phi)=\overline{F}_{\alpha(\Phi)}.
$$
\end{remark}

Let  $\Phi=(\Phi_1|\Phi_2|\cdots|\Phi_k)$ be a set composition, and $\Pi=\Pi_1/\Pi_2/\cdots/\Pi_l$ be a set partition. In the following table, we summarize our results for this section using the operators 
$$\begin{tikzpicture} 
\node at (0,0){$\biguplus,$}; 
\node at (1,0){ $\arx$};
\node at (1.3,-0.15){,};
\node at (2,0){ $\ary$};
\node at (2.3,-0.15){,};
\node at (3,0){ $\arz$};
\end{tikzpicture} $$
on the labelled edge-coloured digraphs in Table \ref{lecd}.

\begin{table}
\caption{Bases for $\NCSym(x)$ and $\NCQSym(x)$ reinterpreted}
\label{table:basesNCQS}
\centering
\begin{tabsize}
 \begin{tabular}{|c|c|c|}
\toprule
{Basis}& {Expression} & {Notation}\\
\midrule
 \begin{tikzpicture} \node at (0,0){Monomial basis of $\NCSym(\bx)$};\end{tikzpicture}&
$ \GY_{ \begin{tikzpicture} 
\node at (-0.6,0){${\Big (}$};
\node at (0,0){$\arx_{i=1}^l$};
\node at (0.5,0){$~~~~C_{\Pi_i}$};
\node at (1.1,0){$\Big )$};
\end{tikzpicture}}(\bx)$  &
 \begin{tikzpicture} \node at (0,0){${\bf m}_\Pi$};\end{tikzpicture} \\
\midrule
\text{Elementary basis of $\NCSym(\bx)$} & $\sum_{(\sigma_1,\ldots,\sigma_l)\in \mathfrak{S}_{\Pi_1}\times \cdots \times \mathfrak{S}_{\Pi_l}} \GY_{(\uplus_{i=1}^l \sigma_i\circ P_{\Pi_i})}(\bx)$&  ${\bf e}_\Pi$\\
\midrule
  \text{Elementary basis of $\NCSym(\bx)$}   & $\GY_{(\biguplus_{i=1}^l K_{\Pi_i})}(\bx)$&${\bf e}_\Pi$\\
\midrule
\text{Complete homogeneous basis of $\NCSym(\bx)$} &$\sum_{(\sigma_1,\ldots,\sigma_l)\in \mathfrak{S}_{\Pi_1}\times \cdots \times \mathfrak{S}_{\Pi_l}} \GY_{(\uplus_{i=1}^l \sigma_i\circ Q_{\Pi_i})}(\bx)$& ${\bf h}_{\Pi}$\\
\midrule
  \text{Power sum basis of $\NCSym(\bx)$} & $\GY_{(\biguplus_{i=1}^l C_{\Pi_i})}(\bx)$& ${\bf p}_\Pi$\\
\midrule
\text{Monomial basis of $\NCQSym(\bx)$}  &$\GY_{ \begin{tikzpicture} 
\node at (-0.6,0){${\Big (}$};
\node at (0,0){$\ary_{i=1}^k$};
\node at (0.5,0){$~~~~C_{\Phi_i}$};
\node at (1.1,0){$\Big )$};
\end{tikzpicture}}(\bx)$ & ${\bf M}_\Phi$\\
\midrule
\text{Fundamental basis of $\NCQSym(\bx)$} &$\GY_{ \begin{tikzpicture} 
\node at (-0.6,0){${\Big (}$};
\node at (0,0){$\ary_{i=1}^k$};
\node at (0.5,0){$~~~~~~Q_{\Phi_i}$};
\node at (1.1,0){$\Big )$};
\end{tikzpicture}}(\bx)$ & ${\bf F}_\Phi$\\
\midrule
\text{Upper-fundamental basis of $\NCQSym(\bx)$} &$\GY_{ \begin{tikzpicture} 
\node at (-0.6,0){${\Big (}$};
\node at (0,0){$\arz_{i=1}^k$};
\node at (0.5,0){$~~~~\fa{C_{\Phi_i}}$};
\node at (1.1,0){$\Big )$};
\end{tikzpicture}}(\bx)$  & $\overline{\bf F}_\Phi$ \\
\bottomrule
\end{tabular}
\end{tabsize}
\end{table}

\fa{
\begin{remark}
The Schur function in noncommuting variables $s_\pi$ in \cite{ALvW} cannot be written as a single  generalized chromatic function of a labelled edge-coloured digraph. For example,  the Schur function in noncommuting variables  \begin{align*}s_{12/34}=&\frac{1}{12} {\bf m}_{1/2/3/4} + \frac{1}{3}{\bf m}_{1/2/34} - \frac{1}{12}{\bf m}_{1/23/4} + \frac{1}{6}{\bf m}_{{1/2 34}} + \frac{1}{12}{\bf m}_{{1/2 4/3}} \\ &+ \frac{1}{6}{\bf m}_{{12}/ {3}/ {4}} + \frac{2}{3}{\bf m}_{{1 2}/ {34}} - \frac{1}{2}{\bf m}_{{1 2 3}/ {4}} + \frac{1}{6}{\bf m}_{{1 2 4}/ {3}} - \frac{1}{12}{\bf m}_{{1 3}/ {2}/ {4}} \\ &- \frac{1}{12}{\bf m}_{{1 3}/ {2 4}} + \frac{1}{6}{\bf m}_{{1 3 4}/ {2}} + \frac{1}{12}{\bf m}_{{1 4}/ {2}/ {3}} - \frac{1}{12}{\bf m}_{{1 4}/ {2 3}},\end{align*} and so in the expansion of $s_{12/34}$ in terms of monomials the coefficient of $\bx_1\bx_1\bx_1\bx_2$ is negative, but the generalized chromatic function of any labelled edge-coloured digraph is a nonnegative linear combination of monomials. 
Nevertheless, the source Schur functions in \cite{ALvW} can be realized as the noncommutative determinant of a matrix where its entries are generalized chromatic functions in noncommuting variables. 
\end{remark}
}

\section{Fundamental Basis to Fundamental Basis: an Injection of the Malvenuto-Reutenauer Hopf Algebra into $\NCQSym(\bx)$}\label{MR}

\fa{Given a list  $i_1i_2\cdots i_n$ of positive integers, define ${\rm std}(i_1i_2\cdots i_n)=u_1u_2\cdots u_n$ where $$u_k=|\{j\in [n]: i_j<i_k\}|+|\{j\in [n]: j\leq k, i_j=i_k\}|,$$ that is  the number of integers in $i_1i_2\cdots i_n$ that are smaller than $i_k$ plus the number of integers equal to $i_k$ that are weakly to the left of $i_k$.}  For example, 
$${\rm std}(3544231)=3756241.$$ Given a permutation $\sigma\in \mathfrak{S}_n$, let $|\sigma|=n.$ Define 
$$F_\sigma=\sum \bx_{i_1}\bx_{i_2}\cdots \bx_{i_n}\in \mathbb{Q}\langle \langle \bx_1,\bx_2,\ldots \rangle \rangle$$  where the sum runs over all lists $i_1i_2\cdots i_n$  with $i_j\in \mathbb{P}$ such that ${\rm std}(i_1i_2\cdots i_n)=\sigma^{-1}$.

The \emph{Malvenuto-Reutenauer} Hopf algebra is the graded Hopf algebra 
$$\displaystyle\SSym(\bx)=\bigoplus_{n\geq 0}\SSym_n(\bx),$$ where $$\SSym_n(\bx)=\mathbb{Q}\text{-span}\{F_{\sigma}:\sigma\in \mathfrak{S}_n\}.$$For more details about the Malvenuto-Reutenauer Hopf algebra, see \cite[Section 8]{GR}.

Given any permutation $\sigma$, let ${\rm Des}(\sigma)=\{i: \sigma(i) > \sigma(i+1)\}$. Let $i_1<i_2<\dots< i_t$ be all elements of ${\rm Des}(\sigma)$. By definition we have that $F_\sigma=\GY_{\bG}(\bx)$ where
 $$
\begin{tikzpicture} 
\node (a) at (-1.9,0){${\bG}=$};
\node (b) at (0,0){$Q_{\{\sigma(1), \dots, \sigma(i_1)\}}$};
\node (c) at (1.5,0){$\ary$};
\node (d) at (3.3,0){$Q_{\{\sigma(i_1+1), \dots, \sigma(i_2)\}}$};
\node (e) at (5,0){$\ary$};
\node (f) at (5.7,0){$\dots$};
\node (g) at (6.3,0){$\ary$};
\node (h) at (8.1,0){$Q_{\{\sigma(i_t+1), \dots, \sigma(n)\}}.$};
\end{tikzpicture} 
$$ 
For the permutation $\sigma$, define the set composition
$$\Phi_\sigma=(\sigma(1) \cdots \sigma(i_1)| \fathree{\sigma(i_1+1)} \cdots \sigma(i_2)| \cdots| \sigma(i_{t-1}+1) \cdots \sigma(i_t)).$$ 
For example, if $\sigma=836791524$, then $${\rm Des}(\sigma)=\{1,5,7\},$$ and 
$F_\sigma=\GY_{\bG}$ where 
$$
\begin{tikzpicture} 
\node (a) at (-1.2,0){${\bG}=$};
\node (b) at (0,0){$Q_{\{8\}}$};
\node (c) at (0.8,0){$\ary$};
\node (d) at (1.9,0){$Q_{\{3,6,7,9\}}$};
\node (e) at (3,0){$\ary$};
\node (f) at (3.9,0){$Q_{\{1,5\}}$};
\node (g) at (4.8,0){$\ary$};
\node (h) at (5.7,0){$Q_{\{2,4\}}.$};
\end{tikzpicture} 
$$ 
Moreover, 
$$\Phi_{\sigma}=(8|3679|15|24).$$

Considering the product and coproduct formulas for the fundamental bases of $\SSym(\bx)$ and $\NCQSym(\bx)$, we have the following injection,
$$
\begin{array}{ccc}
\SSym (\bx)& \rightarrow & \NCQSym(\bx)\\
F_\sigma & \mapsto & F_{\Phi_\sigma}.
\end{array} 
$$

\section{$\NCQSym^r(\bx)$ and Its Bases}\label{sec:basesR}

To conclude we prove that $\NCQSym^r(\bx)$ is a Hopf algebra, where
$${\rm NCQSym}(\bx)={\rm NCQSym}^1(\bx)  \supset {\rm NCQSym}^2(\bx)  \supset \cdots \supset {\rm NCQSym}^\infty(\bx)={\rm NCSym}(\bx).$$ We then establish two natural bases for this Hopf algebra using generalized chromatic functions.

The embedding of $\NCQSym^r(\bx)$ to $\NCQSym(\bx)$ is given by
$${\bf M}_{(\Phi,\Pi)}=\sum_{\Psi\in\Phi\overline{\scriptstyle\dsqcup}\Pi}{\bf M}_{\Psi}$$
where $\scriptstyle\overline{\dsqcup}$ is defined as follows. If $\Phi=(\Phi_1|\Phi_2|\cdots|\Phi_k)$ and $\Pi=\Pi_1/ \Pi_2/\cdots/\Pi_l$, then $\Phi\overline{\scriptstyle\dsqcup}\Pi$ is the set of set compositions $(\Psi_1|\Psi_2|\cdots|\Psi_{k+l})$ such that
$$\{\Psi_1,\Psi_2,\dots,\Psi_{k+l}\}=\{\Phi_1,\Phi_2,\ldots,\Phi_k,\Pi_1,\Pi_2, \ldots,\Pi_l\}$$
and if $\Phi_i=\Psi_p$ and $\Phi_{i+1}=\Psi_q$, then we always have $p<q$.

In the following theorem, we show that $\NCQSym^r(\bx)$ is a Hopf algebra by finding the product and coproduct formulas for the $r$-dominant monomial basis of $\NCQSym ^{r}(\bx)$. 

\begin{theorem}\label{ncqsymr}
For any positive integer $r$, $\NCQSym^r(\bx)$ is a Hopf algebra.
\end{theorem}

\begin{proof}	
\fa{Since $\NCQSym^r(\bx)$ is graded with $\NCQSym^r_0(\bx)=\mathbb{Q}$, by Takeuchi's formula \cite{T},} it is enough to show that $\NCQSym^r(\bx)$ is closed under the product and coproduct of $\NCQSym(\bx)$. 
In $\NCSym(\bx)$, the coproduct on the monomial basis is taking subsets of blocks, followed by standardization.
That is, given a set partition $\Pi=\Pi_1/\Pi_2/\cdots/\Pi_l$, then 
$$\Delta({\bf m}_\Pi)=\sum_{\{i_1,i_2,\ldots, i_l\}\fathree{=} [l]  \atop i_1<i_2<\cdots<i_{j}; i_{j+1}<i_{j+2}<\cdots<i_{l}} {\bf m}_{{\rm std}(\Pi_{i_1}/\Pi_{i_2}/\cdots/\Pi_{i_j})} \otimes  {\bf m}_{{\rm std}(\Pi_{i_{j+1}}/\Pi_{i_{j+2}}/\cdots/\Pi_{i_{l}} )}.$$
Recall that for a set composition $\Phi=(\Phi_1|\Phi_2|\cdots|\Phi_k)$, in $\NCQSym(\bx)$, we have
$$\Delta({\bf M}_\Phi)=\sum_{t=0}^k {\bf M}_{{\rm std}(\Phi_1|\Phi_2|\cdots|\Phi_t)}\otimes {\bf M}_{{\rm std}(\Phi_{t+1}|\Phi_{t+2}|\cdots|\Phi_k)}.$$
Therefore, in $\NCQSym^r(\bx)$, if $(\Phi,\Pi)=(\Phi_1|\Phi_2|\cdots|\Phi_k,\Pi_1/\Pi_2/\cdots/\Pi_l)$, we have
$$\Delta({\bf M}_{(\Phi,\Pi)})=\sum_{\fa{t=0}}^k\sum_{\{i_1,\dots,i_l\}=[l]\atop i_1<\cdots<i_j; i_{j+1}<\cdots<i_l} {\bf M}_{\text{std}(\Phi_{1}|\cdots|\Phi_{t},\Pi_{i_1}/\cdots/\Pi_{i_j})}\otimes {\bf M}_{\text{std}(\Phi_{t+1}|\cdots|\Phi_{k},\Pi_{i_{j+1}}/\cdots/\Pi_{i_l})}.$$
For example, when $r=2$,
$$
\Delta({\bf M}_{((24),1/3)})=1\otimes {\bf M}_{((24),1/3)}+{\bf M}_{(\emptyset,1)}\otimes {\bf M}_{((13),2)}+{\bf M}_{(\emptyset,1)}\otimes {\bf M
}_{((23),1)}+{\bf M}_{((23),1)}\otimes {\bf M}_{(\emptyset,1)}$$
$$
+{\bf M}_{((13),2)}\otimes {\bf M}_{(\emptyset,1)}+{\bf M}_{((12),\emptyset)}\otimes {\bf M}_{(\emptyset,1/2)}+{\bf M}_{(\emptyset,1/2)}\otimes {\bf M}_{((12),\emptyset)}+{\bf M}_{((24),1/3)}\otimes 1.
$$

The product formula is a little more complicated. Before we continue we need to recall the product for the monomial basis of $\NCQSym(\bx)$. Let $\Phi=(\Phi_1|\Phi_2|\cdots|\Phi_k)$ be a set composition of $[n]$, and let $A$ be a subset of $[n]$. \fa{The restriction of $\Phi$ to $A$, $\Phi|_A$, is the set composition obtained by dropping the empty parts of $(\Phi_1\cap A|\Phi_2\cap A|\cdots|\Phi_k\cap A)$.} For example, $$(357|26|14)|_{\{1,3,4\}}=(3|14).$$
Let $\Phi\vDash [n]$ and $\Psi\vDash[m]$.  The {shifted quasi-shuffle} of $\Phi$ and $\Psi$, denoted $\Phi \overrightarrow{\scriptstyle\dsqcup}\Psi$, is the set of set compositions \fa{$\Gamma\vDash [n+m]$} such that  $\Gamma|_{\{1,2,\ldots,n\}}=\Phi$ and ${\rm std}(\Gamma|_{ \{ n+1,n+2,\ldots,n+m\} })=\Psi$. The product formula for the monomial basis of $\NCQSym(\bx)$ is
$${\bf M}_{\Phi}\cdot {\bf M}_{\Psi}=\sum_{\Gamma\in\Phi \overrightarrow{\scriptstyle\dsqcup}\Psi}{\bf M}_{\Gamma}.$$
Now for an $r$-set-composition $(\Phi,\Pi)$, we have
$${\bf M}_{(\Phi,\Pi)}=\sum_{\Psi\in\Phi\overline{\scriptstyle\dsqcup}\Pi}{\bf M}_{\Psi}.$$
Conversely, given a set composition $\Psi$, there is a unique $r$-set-composition $(\Phi,\Pi)$ such that $\Psi\in\Phi\overline{\scriptstyle\dsqcup}\Pi$. We denote $\Phi$ by $\Psi|_{r\comp}$ and $\Pi$ by $\Psi|_{r\pa}$.

Also, for set compositions $\Gamma$ and $\Theta$ we have 
$${\bf M}_{\Gamma}\cdot {\bf M}_{\Theta}=\sum_{\Upsilon\in\Gamma \overrightarrow{\scriptstyle\dsqcup}\Theta}{\bf M}_{\Upsilon}.$$
Conversely, given a set composition $\Upsilon$ with $|\Upsilon|=n+m$, there is a unique pair of set compositions $\Gamma,\Theta$ such that $|\Gamma|=n$, $|\Theta|=m$ and $\Upsilon\in\Gamma\overrightarrow{\scriptstyle\dsqcup}\Theta$. We denote $\Gamma$ by $\Upsilon|_{\{1,\dots,n\}}$ and $\Theta$ by $\Upsilon|_{\{n+1,\dots,n+m\}}$.

Therefore,
$${\bf M}_{(\Phi,\Pi)}\cdot {\bf M}_{(\Psi,\Omega)}=\sum_{\Gamma \in\Phi \overline{\scriptstyle\dsqcup}\Pi, \Theta \in\Psi\overline{\scriptstyle\dsqcup}\Omega}\left(\sum_{\Upsilon\in \Gamma\overrightarrow{\scriptstyle\dsqcup}\Theta}{\bf M}_{\Upsilon}\right)=\sum_{\Upsilon}C_{(\Phi,\Pi),(\Psi,\Omega)}^\Upsilon\bf M_\Upsilon.$$

Let $|(\Phi,\Pi)|=n$ and $|(\Psi,\Omega)|=m$. We first note that all coefficients $C_{(\Phi,\Pi),(\Psi,\Omega)}^\Upsilon$ are $0$ or $1$. Indeed, $C_{(\Phi,\Pi),(\Psi,\Omega)}^\Upsilon=1$ if and only if:
\begin{itemize}
	\item[(i)] $\Phi=(\Upsilon|_{\{1,\dots,n\}})_{r\comp}$.
	\item[(ii)] $\Pi=(\Upsilon|_{\{1,\dots,n\}})_{r\pa}$.
	\item[(iii)] $\Psi=(\Upsilon|_{\{n+1,\dots,n+m\}})_{r\comp}$.
	\item[(iv)] $\Omega=(\Upsilon|_{\{n+1,\dots,n+m\}})_{r\pa}$.
\end{itemize}

Let $C_{(\Phi,\Pi),(\Psi,\Omega)}^{\Upsilon'}=1$ for some $\Upsilon'$, and let $\Phi'=\Upsilon'_{r\comp}$ and $\Pi'=\Upsilon'_{r\pa}$. We want to show that for any $\Upsilon''\in \Phi'\overline{\scriptstyle\dsqcup}\Pi'$, we have $C_{(\Phi,\Pi),(\Psi,\Omega)}^{\Upsilon''}=1$. Then the product ${\bf M}_{(\Phi,\Pi)}\cdot {\bf M}_{(\Psi,\Omega)}$ is in $\NCQSym^r(x)$.

Let $\Gamma'=\Upsilon'|_{\{1,\dots,n\}}$, $\Theta'=\Upsilon'|_{\{n+1,\dots,n+m\}}$, $\Gamma''=\Upsilon''|_{\{1,\dots,n\}}$ and $\Theta''=\Upsilon''|_{\{n+1,\dots,n+m\}}$. Note the fact that in a quasi-shuffle, the blocks of size less than $r$ can only be obtained from blocks of size less than $r$. Since $\Upsilon'_{r\comp}=\Upsilon''_{r\comp}=\Phi'$, we must have $\Gamma'_{r\comp}=\Gamma''_{r\comp}=\Phi$ and $\Theta'_{r\comp}=\Theta''_{r\comp}=\Psi$. And therefore, we must also have $\Gamma'_{r\pa}=\Gamma''_{r\pa}=\Pi$ and $\Theta'_{r\pa}=\Theta''_{r\pa}=\Omega$. Hence, $\Gamma''\in\Phi\overline{\scriptstyle\dsqcup}\Pi$ and $\Theta''\in\Psi\overline{\scriptstyle\dsqcup}\Omega$, that is, $C_{(\Phi,\Pi),(\Psi,\Omega)}^{\Upsilon''}=1$.
\end{proof}

\begin{sloppypar}
Finally we establish two bases for ${\NCQSym^r(\bx)}$. 
Given an $r$-set-composition $(\Phi,\Pi)= ( (\Phi_1|\Phi_2|\cdots|\Phi_k),\Pi_1/\Pi_2/\cdots/\Pi_l)$ of $[n]$, by definition we have that
$${\bf M}_{(\Phi,\Pi)}=
\GY_{
\begin{tikzpicture} 
\node at (-0.6,0){${\Big (}$};
\node at (0,0){$\ary_{i=1}^k$};
\node at (0.5,0){$~~~~C_{\Phi_i}$};
\node at (1.1,0){$\Big )$};
\node at (1.5,0){$\arx$};
\node at (1.9,0){$\Big ($};
\node at (2.5,0){$\arx_{j=1}^l$};
\node at (3,0){$~~~~C_{\Pi_j}$};
\node at (3.6,0){$\Big )$};
\end{tikzpicture}}(\bx).$$
Define 
$$\overline{\bf F}_{(\Phi,\Pi)}=\GY_{
\begin{tikzpicture} 
\node at (-0.6,0){${\Big (}$};
\node at (0,0){$\arz_{i=1}^k$};
\node at (0.5,0){$~~~~C_{\Phi_i}$};
\node at (1.1,0){$\Big )$};
\node at (1.5,0){$\arx$};
\node at (1.9,0){$\Big ($};
\node at (2.5,0){$\arx_{j=1}^l$};
\node at (3,0){$~~~~C_{\Pi_j}$};
\node at (3.6,0){$\Big )$};
\end{tikzpicture}}(\bx).$$
\end{sloppypar}

\begin{proposition}\label{bncqsymr}
Each of the following is a basis for $\NCQSym_n^r(\bx)$.
\begin{itemize}
\item[(i)] $\{ \bM_{(\Phi,\Pi)}: (\Phi,\Pi) \text{~is an $r$-set-composition of $[n]$} \}$.
\item[(ii)] $\{ \overline{\bF}_{(\Phi,\Pi)}:  (\Phi,\Pi) \text{~is an $r$-set-composition of $[n]$} \}$.
\end{itemize} 
\end{proposition}

\begin{proof}
The first set is the $r$-dominant monomial basis for $\NCQSym^r_n(\bx)$. Now, consider that 
$$
\overline{\bF}_{(\Phi,\Pi)}=\sum_{\Psi \text{\svwt{~corrupts~}} \Phi} \bM_{(\Psi,\Pi)}.
$$
Therefore, the second set \fa{is also} a basis for $\NCQSym_n^r(\bx)$.
\end{proof}

\section*{Funding} This work was supported by the National Sciences and Engineering Research Council of Canada. 

\acks{The authors would like to thank Victor Wang for his helpful suggestions.  \faa{They also would like to thank the referees for their insightful recommendations and thoughtful
comments. In particular, they would like to thank one of the
referees for their suggestion on the statement of Proposition \ref{prop:CSasGCF}, and also extend their gratitude to the other referee for suggesting the notation for labelled edge-coloured digraphs and the $\odot$ product, as well as for providing valuable suggestions that helped to simplify the proof of Theorem \ref{thm:gcp}. The use of the $\odot$ product notation played a crucial role in deriving the results presented in Theorem \ref{thm:r-bases}.}
}

\vspace{2cm}

\noindent Farid Aliniaeifard, Department of Mathematics,
 University of British Columbia,
 Vancouver BC V6T 1Z2, Canada, Email: \href{mailto:farid@math.ubc.ca}{farid@math.ubc.ca} \\
 \\
Shu Xiao Li, School of Mathematical Sciences,
Dalian University of Technology,
Dalian Liaoning 116024, P.R. China, Email: \href{mailto:lishuxiao@dlut.edu.cn}{lishuxiao@dlut.edu.cn} \\
\\
Stephanie van Willigenburg, Department of Mathematics,
 University of British Columbia,
 Vancouver BC V6T 1Z2, Canada, Email: \href{mailto:steph@math.ubc.ca}{steph@math.ubc.ca}

\end{document}